\documentclass[11pt]{amsart}
\usepackage[letterpaper, total={6.5in,9in}, includefoot, centering]{geometry}
\usepackage{amsfonts,mathtools}
\usepackage{amsmath,amsthm,amssymb,amscd}
\usepackage[foot]{amsaddr}
\usepackage[T1]{fontenc}
\usepackage{esint}

\usepackage[mathscr]{eucal}
\usepackage{latexsym}
\usepackage[new]{old-arrows}
\usepackage{graphics}
\usepackage{verbatim}
\usepackage{upgreek}
\usepackage{overpic}
\usepackage[all,cmtip]{xy} 
\usepackage{enumitem}
\usepackage{tikz}
\usepackage{tikz-cd}
\usepackage[pagebackref]{hyperref}
\renewcommand*\backref[1]{\ifx#1\relax \else (Page #1) \fi}
\usepackage[font=scriptsize]{caption}
\usepackage{upgreek}
\makeatletter
\newsavebox{\@brx}
\newcommand{\llangle}[1][]{\savebox{\@brx}{\(\m@th{#1\langle}\)}%
  \mathopen{\copy\@brx\kern-0.5\wd\@brx\usebox{\@brx}}}
\newcommand{\rrangle}[1][]{\savebox{\@brx}{\(\m@th{#1\rangle}\)}%
  \mathclose{\copy\@brx\kern-0.5\wd\@brx\usebox{\@brx}}}
\makeatother
\graphicspath{ {./images/} }
\usepackage{tikz}
\usepackage{tikz-cd}
\usepackage{bbm}

\usepackage[authoryear,sort&compress]{natbib}
\usepackage{adjustbox}
\usepackage{float}
\usepackage{algorithm}
\usepackage{algpseudocode}
\usepackage{url}
\bibpunct{[}{]}{;}{n}{,}{,}
\usepackage[caption=false]{subfig}
\usepackage{stackengine}
\usepackage{svg}
\usepackage{cleveref}
\hypersetup{
	colorlinks=true,
    linkcolor=red,
    citecolor=blue,
}

\makeatletter \@addtoreset{equation}{section}

\DeclareMathOperator*{\argmin}{arg\,min}
\newcommand{\prox}{\textnormal{prox}}

\makeatother

\newtheorem{definition}{Definition}[section]
\newtheorem{remark}{Remark}[section]
\newtheorem{example}{Example}[section]

\newtheorem*{prop*}{Proposition}

\newcommand{\bx}{\mathbf{x}}
\newcommand{\bh}{\mathbf{h}}
\newcommand{\by}{\mathbf{y}}

\newcommand{\bv}{\mathbf{v}}

\newcommand{\br}{\mathbf{r}}
\newcommand{\boldf}{\mathbf{f}}

\newcommand{\bsO}{{\boldsymbol{\Omega}}}

\newcommand{\mF}{{\widetilde{F}}}
\newcommand{\wtnabla}{{\widetilde{\nabla}}}

\begin{document}
\title[Nonlinear Splitting for Gradient-Based Unconstrained and Adjoint Optimization]{Nonlinear Splitting for Gradient-Based\\Unconstrained and Adjoint Optimization}

\begin{abstract}
High dimensional and/or nonconvex optimization remains a challenging and important problem across a wide range of fields, such as machine learning, data assimilation, and partial differential equation (PDE) constrained optimization. Here we consider gradient-based methods for solving unconstrained and constrained optimization problems, and introduce the concept of \emph{nonlinear splitting} to improve accuracy and efficiency. For unconstrained optimization, we consider splittings of the gradient to depend on two arguments, leading to semi-implicit gradient optimization algorithms. In the context of adjoint-based constrained optimization, we propose a splitting of the constraint $F(\bx,\theta)$, effectively expanding the space on which we can evaluate the ``gradient''. In both cases, the formalism further allows natural coupling of nonlinearly split optimization methods with acceleration techniques, such as Nesterov or Anderson acceleration. The framework is demonstrated to outperform existing methods in terms of accuracy and/or runtime on a handful of diverse optimization problems. This includes low-dimensional analytic nonconvex functions, high-dimensional nonlinear least squares in quantum tomography, and PDE-constrained optimization of kinetic equations, where the total number of high-dimensional kinetic solves is reduced by a factor of three compared with standard adjoint optimization.  
\end{abstract}

\author{Brian K. Tran$^{1,5}$, Ben S. Southworth$^1$, David B. Cavender$^2$, \\ Sam Olivier$^3$, Syed A. Shah$^1$, Tommaso Buvoli$^4$}
\address{$^1$Los Alamos National Laboratory, Theoretical Division, Los Alamos, New Mexico 87545}
\address{$^2$University of California San Diego, Department of Mathematics, La Jolla, California 92093}
\address{$^3$Los Alamos National Laboratory, Computer, Computational and Statistical Sciences Division, Los Alamos, New Mexico 87545}
\address{$^4$Tulane University, Department of Mathematics, New Orleans, LA 70118}
\address{$^5$University of Colorado Boulder, Department of Applied Mathematics, Boulder, CO 80303}
\address{\emph{E-mail}: \protect{brian.tran@colorado.edu}, \protect{southworth@lanl.gov}, \protect{dcavende@ucsd.edu}} 
\address{\emph{E-mail}: \protect{solivier@lanl.gov}, \protect{shah@lanl.gov}, \protect{tbuvoli@tulane.edu}  }
\allowdisplaybreaks

\maketitle
{\hypersetup{linkcolor=black}\tableofcontents}

\section{Introduction}\label{sec:intro}

High dimensional and/or nonconvex optimization remains a challenging and important problem across a wide range of fields, such as machine learning, data assimilation, and partial differential equation (PDE) constrained optimization. This paper is concerned with gradient-based methods for solving two classical optimization problems, namely the unconstrained optimization problem 
\begin{align}
	\min_{\bx} C(\bx),	
	\label{eq:opt-problem-unconstrained}
\end{align}
and the constrained optimization problem
\begin{align}
    &\min_{\theta} J(\bx;\theta), \quad
    \text{subject to }  F(\bx;\theta) = 0	.
\label{eq:opt-problem-constrained}
\end{align}%
For both problems, we will employ {\em nonlinear splittings} to improve the efficiency and accuracy of gradient-based optimization. Specifically, given a function $f(\bx)$, we define a nonlinear splitting of $f$ to be a function $g(\bx,\by)$ that satisfies $g(\bx,\bx) = f(\bx)$ and is differentiable in both arguments. In the context of unconstrained optimization, we will consider splittings of the gradient $\nabla C$ to depend on two arguments, leading to semi-implicit gradient optimization algorithms. For the constrained optimization \eqref{eq:opt-problem-constrained}, the constrained objective is often sufficiently complex that treating components of the gradient implicitly may not be tractable. Instead, we propose a splitting of the constraint $F(\bx;\theta)$ within an adjoint-based optimization. This effectively expands the space on which we can evaluate the ``gradient,'' and induces a sort of splitting in the gradient as well. In both cases, the formalism further allows natural coupling of nonlinearly split optimization methods with acceleration techniques, such as Nesterov or Anderson acceleration. 


This paper is organized as follows. We begin by motivating nonlinear splitting in the context of gradient flow and classical algorithms like gradient descent, proximal point method, and proximal gradient method in \Cref{sec:background}. In \Cref{sec:main}, we develop the nonlinear splitting framework for gradient-based optimization. \Cref{sec:unconstrained} focuses on the unconstrained setting, where we introduce the notion of a nonlinear splitting of a gradient and use this to define the nonlinear splitting gradient descent method. In \Cref{sec:acceleration-unconstrained}, we discuss momentum-based and fixed-point-based acceleration of nonlinear splitting gradient descent. Throughout, we discuss connections to existing methods and generalizations thereof. In \Cref{sec:equality-constraints}, we discuss nonlinear splitting in the setting of equality-constrained optimization. Here, we introduce the notion of a nonlinear splitting of the equality constraint and a nonlinear splitting of the gradient. Using a modified adjoint approach, we show how to systematically construct a nonlinear splitting of the gradient, given a nonlinear splitting of the equality constraint. A general framework is provided where any nonlinear approximation or solver written in fixed-point form can naturally induce such a splitting. In \Cref{sec:acceleration}, we provide momentum-based and fixed-point-based accelerated variants of the equality-constrained nonlinear splitting gradient method. 

In \Cref{sec:numerical-example}, we demonstrate these proposed methods on a variety of numerical examples. For unconstrained optimization, we consider two examples in \Cref{sec:uncons-examples} that highlight the advantage of the linear stability gained by using nonlinear splitting gradient descent. Specifically, in \Cref{sec:benchmarks}, we consider several two-dimensional non-convex benchmark unconstrained optimization problems, and demonstrate how the linear stability from these nonlinear splitting methods allows one to take large steps in order to escape local minima and find the global minima of these non-convex functions. In terms of objective value versus runtime, these nonlinear splitting methods outperform most of MATLAB's built-in optimizers, with MATLAB's global search algorithm being the only one that is comparable. That being said, the nonlinear splitting methods are implemented via naive looping in MATLAB scripts and thus we expect some performance degradation compared with the optimized methods provided in MATLAB. In \Cref{sec:NLS}, we consider a nonlinear least squares problem to reconstruct a pure quantum state arising in quantum tomography. This example demonstrates linearization of the nonlinear splitting methods, using a single Newton iteration, as well as Anderson acceleration of the corresponding method. For equality-constrained optimization, we consider two examples in \Cref{sec:cons-examples}. In \Cref{sec:neutron-transport}, we consider a source optimization problem for the nonlinear neutron transport equation, where we demonstrate that our proposed methods significantly reduce the number of linear transport solves in order to achieve the same objective value and equality constraint residual, compared to their standard counterparts. In \Cref{sec:time-int-example}, we show how the nonlinear splitting framework can be applied in the case when the equality constraint is given by the dynamics of an ordinary differential equation (ODE), which leads us to a momentum-based scheme for problems with time dependence.

\section{Gradient flow and nonlinear splitting}\label{sec:background}

The idea of using ordinary differential equation (ODE) theory and techniques to solve nonlinear root-finding and optimization problems for cost function $C(\bx)$ dates back many decades, and many optimization methods can be seen as certain discrete approximations to continuous gradient flow equations,
\begin{equation}\label{eq:gradient-flow}
    \dot{\bx} = -\nabla C(\bx),
\end{equation}
or continuous modifications thereof. It is also common for optimization algorithms to be based on solving local optimization problems during each step/iteration towards the global solution. For example, suppose $C(\bx)$ has a Lipschitz continuous differential $DC(\bx)$ with constant $L$, then \cite[Lemma 54]{Kelner.2013} $C(\bx+\bh) \leq C(\bx) + DC(\bx)\bh + \frac{L}{2}\|\bh\|^2$. Taking a greedy minimization with respect to $\bh$ yields the gradient descent algorithm,
\begin{align}\label{eq:GD-standard}
    \bx^{k+1} & = \bx^k + \textnormal{arg}\min_{\bh} \Big( C(\bx) + \langle \nabla C(\bx),\bh\rangle  + \frac{L}{2}\|\bh\|^2 \Big) 
    = \bx^k - \tfrac{1}{L}\nabla C(\bx^k),
\end{align} 
which is exactly a forward Euler discretization of the continuous gradient flow \eqref{eq:gradient-flow}. If we have a matrix-induced $A$-norm for symmetric positive definite $A$ (or natural gradients with a locally defined inner product \cite{amari1998natural}), $\|\bh\|_A^2 = \bh^TA\bh$, and minimization with respect to $\bh$ yields the correction
\begin{equation}\label{eq:lin-prec-GD}
    \bx^{k+1} = \bx^k - \tfrac{1}{L}A^{-1}\nabla C(\bx^k),
\end{equation}
which is a linearly preconditioned gradient descent \cite{precGD1, precGD2, precGD3}. This remains an explicit approximation to the continuous gradient flow, except here we take the gradient with respect to the $A$ inner product, where $\nabla_A C(\bx) = A^{-1}\nabla C(\bx)$. Newton's method can be derived from continuous gradient flow with respect to a Hessian-induced inner product (or duality pairing for non-SPD Hessians \cite{TrSoBlOl2025}), and Levenberg–Marquardt follows from explicit discrete approximations to the modified continuous gradient flow $\dot{\bx} = -H(\bx)^{-1}\nabla B(\bx)$, e.g. \cite{Abbott.1975,Botsaris.1978,D..1998}. 

Damping or momentum can also be introduced via second-order formulations. The heavy-ball with friction method, which has had a resurgence of interest due to the development of large-scale optimization and machine learning problems, was developed in 1964 \cite{Polyak.1964}. There, they consider multistep iterative methods to solve linear and nonlinear systems, and relate the multistep methods to the continuous second-order ODE
\begin{equation}\label{eq:cont-heavy-ball}
    \ddot{\bx} + \alpha_1\dot{\bx} + \alpha_2\nabla C(\bx) = 0,
\end{equation}
for objective function $C(\bx)$. More general second order equations for optimization are typically of the form 
$
    \mu(t)\ddot{\bx} + \beta(t)\dot{\bx} + \nabla C(\bx)/\mu(t) = 0
$, e.g., \cite{Incerti.1979,Aluffi-Pentini.1984,Brown.1989,Alvarez.2000,ATTOUCH.2000}.
The continuous limit of Nesterov's accelerated gradient (NAG) descent, for example, is given by 
$
    \ddot{\bx} + \frac{3}{t}\dot{\bx} + \nabla  C(\bx) = 0,
$
with $\bx(0) = \bx_0, \dot{\bx}(0) = 0$ and a change of variable in time to decay the learning rate \cite{Su.2015}.

Above are all explicit approximations to continuous gradient flows. Although each step can be posed as a certain local optimization problem, the solution is either analytic or linear in nature. 
In contrast, proximal methods are based around optimizing or solving some ``proximal,'' i.e. nearby problem. The typical $\ell^2$ proximal operator $\prox_F:\mathbb{R}^n\to\mathbb{R}^n$ with constant $\gamma_k > 0$ takes the form 
\begin{equation}\label{eq:prox}
    \prox_{\gamma C}(\bx) = \argmin_{\bv}\left( C(\bv) + \tfrac{1}{2\gamma}\|\bx-\bv\|^2\right).
\end{equation}
Much of the literature on proximal methods is specifically for convex functions, where $\bx^*$ minimizes $C(\bx)$ if and only if $\bx^* = \prox_{C}(\bx^*)$, suggesting the fixed-point algorithm \cite{Rockafellar.1976}
\begin{equation}\label{eq:prox-point}
    \bx^{k+1} = \prox_{\gamma_k  C}(\bx^k).
\end{equation}
This is known as the proximal point method. If we assume that $C(\bx)$ is a closed convex function with unique minimum such that $\nabla C(\bx^*) = \mathbf{0}$ (or we can simply relax the objective to finding a local rather than global minimum in \eqref{eq:prox}), note that $\prox_{\gamma_k C}(\bx) = \left\{ \bv\textnormal{ : } \bv + \gamma_k\nabla C(\bv) = \bx \right\}$. Then the proximal point algorithm \eqref{eq:prox-point} corresponds to an \emph{implicit} backward Euler discretization of the continuous gradient flow \eqref{eq:gradient-flow} \cite{Luo.2022}:
\begin{equation}\label{eq:gradient-flow-BE}
        \bx^{k+1} + \gamma_k\nabla C(\bx^{k+1}) = \bx^k.
\end{equation}

Combining these two approaches, the proximal gradient method considers a bi-objective minimization $\min_{\bx} C_1(\bx) + C_2(\bx)$, where $C_1$ and $C_2$ are typically assumed to be convex, and $C_2(\bx)$ is usually a regularization in the unconstrained minimization of $C_1(\bx)$. The proximal gradient method then defines iterations \cite{Parikh.2014}
\begin{equation}\label{eq:prox-gradient1}
    \bx^{k+1} = \prox_{\gamma_k C_2}(\bx^k - \gamma_k\nabla C_1(\bx^k))
        = \argmin_{\bv}\left( C_2(\bv) + \tfrac{1}{2\gamma_k} \|\bx^k - \gamma_k\nabla C_1(\bx^k)-\bv\|^2\right).
\end{equation}
Expressing the $\argmin_{\bv}$ in \eqref{eq:prox-gradient1} as finding a zero gradient (again either by assumption of convexity or relaxing the minimization to finding a local minimum) yields
$
    \bx^{k+1} = \{ \bv \textnormal{ : } \bv + \gamma_k\nabla C_2(\bv) = \bx^k - \gamma_k\nabla C_1(\bx^k) \}.
$
This exactly corresponds to a first-order implicit-explicit formulation (also referred to as a forward-backward splitting) applied to the continuous gradient flow \eqref{eq:gradient-flow} for $C(\bx) = C_1(\bx)+C_2(\bx)$:
\begin{equation}\label{eq:gradient-flow-IMEX}
    \bx^{k+1} + \gamma_k \nabla C_2(\bx^{k+1})= \bx^k - \gamma_k \nabla C_1(\bx^k).
\end{equation}
Note that, when $C_2 \equiv 0$, this reduces to gradient descent for $C_1$ \eqref{eq:GD-standard}, and when $C_1 \equiv 0$, this reduces to the proximal point algorithm for $C_2$ \eqref{eq:prox-point}.

Implicit-explicit and splitting-based time integration schemes have been widely successful for simulation of partial differential equations, particularly when there are multiple coupled scales or physics (e.g., \cite{ascher1995implicit}). Broadly, explicit schemes (like gradient descent) can be intractably slow due to stability restrictions on the stepsize, while implicit schemes (like the proximal point method) are too expensive due to the required nonlinear implicit solve at each step. Implicit-explicit (IMEX) methods attempt to find a balance between the two for efficient (or tractable) simulation. The art of IMEX integration in practice is determining an appropriate splitting (or ``partitioning'' as it is often referred to in the time integration literature) of the problem that ensures stable and accurate integration. The overwhelming majority of literature on splitting in numerical PDEs assumes a linear additive or component splitting of the problem, as in the proximal gradient method \eqref{eq:gradient-flow-IMEX}. However, not all problems have a natural linear splitting, something discussed in detail for thermal radiation transport (TRT) and radiation hydrodynamics in \cite{imex-trt,imex-radhydro}. This motivated the development of a new class of time integration methods for \emph{nonlinear splitting} of operators \cite{nprk1,nprk2}, which follows from the semi-implicit integrators of \cite{Boscarino.2016} and a certain approach used for decades in the TRT literature \cite{Larsen.1988}. For the ODE $\dot{\by} = F(\by)$, one first defines a nonlinear splitting of $F$, which is a function $\mF$ of two arguments such that $\widetilde{F}(\by, \by) = F(\by)$. The nonlinearly split Euler method for approximately evolving the  the system in time is then
\begin{equation}\label{eq:nprk1}
    \by^{k+1} = \by^k + \Delta t \mF(\by^k,\by^{k+1}).
\end{equation}
The timestep \eqref{eq:nprk1} is explicit in the first argument of $\widetilde{F}$ and implicit in the second argument. Unlike the fully implicit method \eqref{eq:gradient-flow-BE}, nonlinear splitting allows one to distribute the implicitness in a more computationally efficient fashion. For example, if $F(y) = y^2$, we can define $\mF(u,v) = uv \implies \mF(y,y) = y^2 = F(y)$. Then each step of \eqref{eq:nprk1} requires a linear solve in $v$, but not a nonlinear solve.

If $\widetilde{F}$ is nonlinear in its second argument, then the timestep \eqref{eq:nprk1} requires a nonlinear solve. By performing only a single Newton iteration as an inexact solve, we obtain the Rosenbrock method  
\begin{equation}\label{eq:nprk2}
    \by^{k+1} = \by^k + \Delta t J_2[\by_k]^{-1} \mF(\by^k,\by^{k}),
\end{equation}
where $J_2[\by_k] = D_2F(\by_k,\by_k)$ denotes the Fr\'{e}chet derivative of $\mF$ with respect to the second argument, evaluated at $\by_k$ in both arguments. Replacing the fully nonlinear solve with a single linear solve does not affect the order of accuracy and both \eqref{eq:nprk1} and \eqref{eq:nprk2} are first-order integrators for approximating the exact flow. Moreover, in \eqref{eq:nprk2}, the nonlinear splitting takes a more subtle role in defining the linear preconditioner $J_2$ that is used to stabilize the otherwise explicit update. In optimization, where it is critical to minimize cost per iteration in method design, only requiring a single linear implicit solve is appealing, and later we demonstrate this to be effective on unconstrained nonlinear least squares problems. The unconstrained optimization methods proposed in \Cref{sec:unconstrained} can largely be seen as generalizations or modifications of \eqref{eq:nprk1}. 

Lastly, we remark that higher order nonlinearly split methods have been developed for the purposes of numerical PDEs \cite{nprk1,nprk2}, however, it is not clear the benefit of higher order methods for gradient flow interpretations of optimization problems.

\section{Gradient-based optimization with nonlinear splitting}\label{sec:main}


\subsection{Nonlinear splitting for unconstrained optimization}\label{sec:unconstrained}
In this section, we consider unconstrained optimization, i.e., to minimize a function $C: \mathbb{R}^N \rightarrow \mathbb{R}$, which we assume to be $C^1$. The gradient descent (GD) method is given by \eqref{eq:GD-standard} with stepsize $1/L\mapsto \gamma_k > 0$, and the method is specified with an initial iterate $\bx^0 \in \mathbb{R}^N$. We now introduce the notion of a nonlinear splitting of the gradient. 

\begin{definition}\label{def:nonlinear-splitting-unconstrained}
    A \textbf{nonlinear splitting of the gradient} of $C: \mathbb{R}^N \rightarrow \mathbb{R}$ (or simply, a nonlinearly split gradient) is a mapping $\widetilde{\nabla}C: \mathbb{R}^N \times \mathbb{R}^N \rightarrow \mathbb{R}^N$ such that
    \begin{equation}\label{eq:gradient-consistency}
       \wtnabla C(\bx,\bx) = \nabla C(\bx)
    \end{equation}
    for all $\bx \in \mathbb{R}^N$. Furthermore, we assume that $\wtnabla C$ is continuously differentiable in both arguments.
\end{definition}
%
%
Given a nonlinearly split gradient $\widetilde{\nabla}C$ of $C$, we then define the nonlinearly split gradient descent (NS-GD method) as
\begin{equation}\label{eq:GD-partitioned}
    \bx^{k+1} = \bx^k - \gamma_k \widetilde{\nabla} C(\bx^{k},\bx^{k+1}).
\end{equation}

As a first example of a nonlinear splitting of a gradient, we define a linearly-implicit nonlinear splitting. The basic idea is to choose a splitting such that the implicit equation \eqref{eq:GD-partitioned} for $\bx^{k+1}$ for NS-GD only requires a linear solve. This is inspired from the time integration literature, particularly, linearly-implicit Runge--Kutta schemes \cite{Calvo.2001, Boscarino.2015}.

\begin{example}[Linearly-implicit nonlinear splitting of a gradient]\label{ex:li-gradient}
    Suppose that the gradient of $C$ can be expressed in the following form
    \begin{equation}\label{eq:grad-C-linear-form}
        \nabla C(\bx) = L(\bx) \bx + g(\bx),
    \end{equation}
    where $L: \mathbb{R}^N \rightarrow \mathbb{R}^{N \times N}$ is a state-dependent linear operator and $g: \mathbb{R}^N \rightarrow \mathbb{R}^N$ is a generally nonlinear map. The linearly-implicit nonlinear splitting of the gradient is given by
    \begin{equation}\label{eq:LI-partition}
        \widetilde{\nabla}C(\bx,\by)  \coloneqq  L(\bx) \by + g(\bx).
    \end{equation}
    We refer to this as a linearly-implicit nonlinear splitting since, in the \textup{NS-GD} \eqref{eq:GD-partitioned} method, a linear solve is needed to obtain the updated iterate. Specifically, with this choice of splitting, the \textup{NS-GD} method \eqref{eq:GD-partitioned} becomes
    $$ \bx^{k+1} = \bx^k - \gamma_k ( L(\bx^k)\bx^{k+1} + g(\bx^k) ), $$
    and thus, the updated iterate is obtained by a single linear solve
    $$ \bx^{k+1} = (I + \gamma_k L(\bx^k))^{-1} (\bx^k - \gamma_k g(\bx^k)). $$
\end{example}

For fixed $\bx$, viewing $\widetilde{\nabla}C(\bx,\by)$ as a one-form with respect to $\by$, if its exterior derivative is zero, then there is an associated cost function $D_\bx(\by)$ such that ${\nabla_{\by} D_\bx(\by) = \widetilde{\nabla}C(\bx,\by)}$. In the case of a linearly-implicit splitting, this implies that there exist matrices $\mathbf{A}(\bx)$, $\mathbf{B}(\bx)$ and vector $\mathbf{b}$ such that
\begin{align}
    \widetilde{\nabla}C(\bx,\by) = \nabla_\by D_\bx(\by) 
    \quad \text{for} \quad   
    D_\bx(\by) = \by^T \mathbf{A}(\bx) \by + \mathbf{B}(\bx) \by + \mathbf{b}.
\end{align}
In this scenario, each step of NS-GD applied to $\widetilde{\nabla}C(\bx,\by)$ is equivalent to solving a new quadratic optimization problem using implicit gradient descent. If the nonlinear splitting is selected so that the  minimum of each quadratic function $D(\by; \bx)$  is also the global minimum of $C(\bx)$ (in general highly non-trivial to do), then NP-GD will escape all local minimum and move towards the global min.

For example, consider the simplified case where a nonlinear split gradient has the form
\begin{align*}
    \wtnabla C(\bx,\by) &= \mathbf{a}(\mathbf{x}) \odot \by,
\end{align*}	
where $\odot$ denotes the Hadamard product. Each iteration of NS-GD is equivalent to implicit gradient descent applied to the function
\begin{align*}
    C(\mathbf{x},\mathbf{y}) = \frac{1}{2} a(\bx)^T (\by \odot \by)
\end{align*}
For all $\bx$ for which $a(\bx)_i > 0$, the iteration will move the solution towards zero. An example of this phenomena can be seen in our numerical experiments on the Rastrigin function in \Cref{sec:benchmarks} where NS-GD obtains the global minimum in less than 10 iterations.

As we previously alluded to in \Cref{sec:intro}, the nonlinear splitting framework defines a broad class of gradient-based descent methods that include existing methods, for example linear additive splitting such as in the proximal gradient method as discussed in \Cref{sec:background}. As another example, the discrete gradient method \cite{discgrad} can be interpreted as a nonlinear splitting gradient descent method.

\begin{example}[Discrete gradient methods]\label{ex:disc-grad}
    Following \cite{discgrad}, a discrete gradient is precisely a nonlinear splitting of the gradient $\wtnabla C(\bx,\by)$ with an additional requirement that the splitting satisfies the mean value property
    \begin{equation}\label{eq:mean-value-property}
        \langle \wtnabla C(\bx,\by), \by - \bx\rangle = C(\by) - C(\bx),
    \end{equation}
    where $\langle\cdot,\cdot\rangle$ denotes the standard inner product on $\mathbb{R}^N$. The additional requirement \eqref{eq:mean-value-property} ensures that the nonlinear splitting gradient descent step \eqref{eq:GD-partitioned} is monotone decreasing in the cost, i.e.,
    \begin{equation}\label{eq:cost-monotonocity}
        C(\bx^{k+1}) \leq C(\bx^{k}),
    \end{equation}
    for any stepsize $\gamma_k \geq 0$. 
    
    An example of a discrete gradient is the average vector field of $\nabla C$ defined by
    \begin{equation}\label{eq:AVF-partition}
        \widetilde{\nabla}C(\bx,\by)  \coloneqq  \int_0^1 \nabla C(\xi \bx + (1-\xi) \by) d\xi.
    \end{equation}
    $\widetilde{\nabla}C(\bx,\by)$ is interpreted as the average of the gradient over the line segment connecting $\bx$ and $\by$.

\end{example}

\subsubsection{Acceleration}\label{sec:acceleration-unconstrained}
We discuss two types of acceleration of the NS-GD method \eqref{eq:GD-partitioned}, by using momentum and by viewing the NS-GD method as a fixed point iteration.

\textbf{Momentum.} The Nesterov accelerated gradient (NAG) method \cite{Nes83, Nes04} is given by
\begin{subequations}\label{eq:NAG-standard}
\begin{align}
    \bx^{k+1} &= \bx^k + \mu \bv^{k} - \gamma_k \nabla C(\bx^k + \mu \bv^{k}), \\
    \bv^{k+1} &= \bx^{k+1} - \bx^k,
\end{align}
\end{subequations}
where $\mu > 0$ is the momentum factor and $\gamma_k>0$ is again the stepsize; the method is initialized with an initial iterate $\bx^0 \in \mathbb{R}^N$ and an initial step direction $\bv^0 = 0 \in \mathbb{R}^N$. The NAG method \eqref{eq:NAG-standard} can be interpreted as first doing a ``look-ahead" step $\bx^k \mapsto \bx^k + \mu \bx^k$, followed by a gradient step, with the gradient evaluated at the look-ahead point, $\nabla C(\bx^k + \mu \bv^{k})$. We express NAG in the form \eqref{eq:NAG-standard} since the generalization to a nonlinearly splitting form is straightforward, although there are other equivalent formulations.

By replacing the gradient appearing in NAG by the nonlinearly split gradient, the nonlinearly split NAG (NS-NAG) method is given by
\begin{subequations}\label{eq:NAG-partitioned}
\begin{align}
    \bx^{k+1} &= \bx^k + \mu \bv^{k} - \gamma_k \widetilde{\nabla} C(\bx_k + \mu \bv_{k}, \bx^{k+1}), \\
    \bv^{k+1} &= \bx^{k+1} - \bx^k.
\end{align}
\end{subequations}

\sloppy \textbf{Fixed-points and linear preconditioning from a nonlinear splitting.} Another approach for accelerating gradient-based descent is to view GD as a fixed-point iteration $\bx^{k+1} = G(\bx^k)$, with fixed-point operator
\begin{equation}\label{eq:GD-fixed-point-operator}
    G(\bx)  \coloneqq  \bx - \gamma \nabla C(\bx).
\end{equation}
A fixed-point of $G$ then corresponds to a stationary point of $C$. One can also generalize and consider a fixed-point operator with linear preconditioning of the gradient
$$ G_A(\bx) := \bx - \gamma A(\bx)^{-1} \nabla C(\bx), $$
where we assume $A(\bx)$ is invertible. Fixed-point iterations of $G_A$ correspond to the linearly preconditioned gradient descent method \eqref{eq:lin-prec-GD}. Note that the fixed-points of $G_A$ are precisely the same as the fixed-points of $G$; however, the convergence properties of the fixed-point iterations are generally different, and $A(\bx)$ is usually chosen to improve convergence and/or stability \cite{precGD1, precGD3}.

As discussed in \Cref{sec:background}, one can consider a linearized solve of a nonlinear splitting method resulting in a linearly-implicit Rosenbrock method \eqref{eq:nprk2}. Here, we will show that such an approach for a NS-GD step can be interpreted as a linear preconditioning of gradient descent. Performing a single Newton iteration, linearized about $\bx^k$, for the nonlinear solve required in the NS-GD step \eqref{eq:GD-partitioned} yields precisely a linearly preconditioned gradient descent method,
\begin{equation}\label{eq:GD-splitting-Newton-1}
    \bx^{k+1} = \widetilde{G}(\bx^k) := \bx^k - \gamma (I + \gamma D_2 \wtnabla C(\bx^k, \bx^{k}))^{-1} \nabla C(\bx^{k}),
\end{equation}
where $\widetilde{G}$ denotes the corresponding fixed-point operator. Note that the preconditioner $A(\bx) := I + \gamma D_2 \wtnabla C(\bx, \bx)$ is defined from the choice of nonlinear splitting. We refer to \eqref{eq:GD-splitting-Newton-1} as NS-GD-Newton(1).

Naturally, one can consider accelerated variants of \eqref{eq:GD-splitting-Newton-1}. For example, using Nesterov momentum, an accelerated method can be constructed by instead linearizing the Newton iteration about the look-ahead point $\overline{\bx}^k  \coloneqq  \bx^k + \mu \bv^k$, which yields
\begin{subequations}\label{eq:NS-NAG-Newton-1}
\begin{align}
    \bx^{k+1} &= \overline{\bx}^k - \gamma (I + \gamma D_2 \wtnabla C(\overline{\bx}^k, \overline{\bx}^k))^{-1} \nabla C(\overline{\bx}^k), \\
    \bv^{k+1} &= \bx^{k+1} - \bx^k.
\end{align}
\end{subequations}
We refer to \eqref{eq:NS-NAG-Newton-1} as NS-NAG-Newton$(1)$. As before, the fixed-points of the iteration \eqref{eq:NS-NAG-Newton-1} are the same as those of $G$ and $\widetilde{G}$. Alternatively, one can consider acceleration methods applied to the nonlinear splitting fixed-point operator $\widetilde{G}$ defined in \eqref{eq:GD-splitting-Newton-1}. For example, a common method for accelerating fixed-point iterations is Anderson Acceleration (AA) \cite{AA1, AA2, Saad_2025}. Let $g(\bx) := \widetilde{G}(\bx) - \bx$, $\bx^1 = \widetilde{G}(\bx^0)$, and $m \geq 1$ be an integer. AA applied to the fixed-point operator $\widetilde{G}$, for $k=1,2,\dots$, is given by
\begin{subequations}\label{eq:NS-AA-Newton-1}
\begin{align}
    m_k &= \min\{m,k\}, \\
    \mathcal{G}^k &= [g(\bx^{k-m_k}) \ \dots\ g(\bx^{k})], \\
    \boldsymbol{\eta}^k &= \argmin_{\boldsymbol{\eta} = (\eta_0,\dots,\eta_{m_k}) : \sum_i \eta_i =1 } \|\mathcal{G}^k\boldsymbol{\eta}\|_{l^2}, \label{eq:AA-ls}\\
    \bx^{k+1} &= \sum_{i=0}^{m_k} \boldsymbol{\eta}^k_i \widetilde{G}(\bx^{k-m_k+i}).
\end{align}
\end{subequations}
We refer to \eqref{eq:NS-AA-Newton-1} as NS-AA-Newton(1). Compared with NS-GD-Newton(1) \eqref{eq:GD-splitting-Newton-1}, NS-AA-Newton(1) has an additional least squares solve per iteration; the least squares solve in \eqref{eq:AA-ls} is assumed to be very low dimensional compared to the parameter space, $m\ll N$, and thus the total cost of the AA update remains $\mathcal{O}(N)$ for $\bx\in\mathbb{R}^N$, a cost that becomes marginal with more complex preconditioners or gradient computations.

We conclude this section with another example, where we show that, for nonlinear least squares with a particular choice of nonlinear splitting, NS-GD-Newton$(1)$ is equivalent to the Levenberg--Marquardt (LM) algorithm \cite{LM1, LM2}.

\begin{example}[Levenberg–Marquardt algorithm as NS-GD-Newton$(1)$]\label{ex:LM-splitting}
    Consider the nonlinear least-squares problem, given by
    \begin{equation}\label{eq:NLS-cost}
       \min_{\bx \in \mathbb{R}^M} C(\bx)  \coloneqq  \frac{1}{2}\sum_{i=1}^N r_i(\bx)^2 = \frac{1}{2}\sum_{i=1}^N (f_i(\bx) - v_i)^2 = \frac{1}{2}\|\boldf(\bx) - \bv\|^2_{l^2(\mathbb{R}^N)},
    \end{equation}
    where $\bv \in \mathbb{R}^N$ is fixed, $\boldf: \mathbb{R}^M \rightarrow \mathbb{R}^N$ is a differentiable and generally nonlinear mapping, and the residual is defined as $\br(\bx)  \coloneqq  \boldf(\bx) - \bv$. Denoting the Jacobian of $\br$ (equivalently, of $\boldf$) as
    $$ \mathbb{J}(\bx)  \coloneqq  \frac{\partial \br(\bx)}{\partial \bx}, $$
    the Levenberg--Marquardt algorithm is given by
    \begin{equation}\label{eq:LM-algorithm}
        (\mathbb{J}^T_k\mathbb{J}_k + \sigma_k I)(\bx^{k+1} - \bx^k) = -\mathbb{J}^T_k \br(x^k),
    \end{equation}
    where $\mathbb{J}_k  \coloneqq  \mathbb{J}(\bx^k)$ denotes the Jacobian of the residual at the $k^{th}$ iterate and $\sigma_k > 0$ is a damping parameter.

    Now, we show that the LM algorithm is equivalent to \textup{NS-GD-Newton}$(1)$  with the following choice of nonlinear splitting, after an appropriate rescaling. The gradient of the nonlinear least squares cost is
    $$ \nabla C(\bx) = \mathbb{J}(\bx)^T \br(\bx). $$
    We define a nonlinear splitting of this gradient by
    \begin{equation}\label{eq:nonlin-LS-split}
        \wtnabla C(\bx, \by) = \mathbb{J}(\bx)^T \br(\by).
    \end{equation}
    Then, $D_2\wtnabla C(\bx, \by) = \mathbb{J}(\bx)^T \mathbb{J}(\by)$. In this case, \eqref{eq:GD-splitting-Newton-1} becomes 
    $$ (I + \gamma_k \mathbb{J}^T_k \mathbb{J}_k ) (\bx^{k+1} - \bx^k) =  - \gamma_k \mathbb{J}^T_k \br(\bx^{k  }) $$
    Dividing both sides by $\gamma_k > 0$ and defining $\sigma_k = 1/\gamma_k$ yields precisely the LM algorithm \eqref{eq:LM-algorithm}.

    If we consider instead \textup{NS-NAG-Newton}$(1)$, we have with this choice of nonlinear splitting
    \begin{subequations}
    \begin{align}
        (\gamma_k\overline{\mathbb{J}}_k^T \overline{\mathbb{J}}_k + I ) (\bx^{k+1} - \overline{\bx}^k) &= - \gamma_k \overline{\mathbb{J}}_k^T r(\overline{\bx}^{k}), \\
        \bv^{k+1} &= \bx^{k+1} - \bx^k,
    \end{align}
    \end{subequations}
    where $\overline{\bx}^k = \bx_k + \mu \bv_{k}$ again denotes the look-ahead point and $\overline{\mathbb{J}}_k = \mathbb{J}(\overline{\bx}^k)$ denotes the Jacobian linearized about the look-ahead point. This reproduces the ``inertial Levenberg--Marquardt" algorithm recently proposed in \cite{inertialLM}. Similarly, we can view LM as a fixed-point algorithm and accelerate it with AA. We will compare \textup{NS-GD-Newton}$(1)$, \textup{NS-NAG-Newton}$(1)$, and \textup{NS-AA-Newton}$(1)$ for a nonlinear least squares problem in the second numerical example, \Cref{sec:NLS}.
\end{example}

\subsection{Nonlinear splitting with equality constraints}\label{sec:equality-constraints}
Throughout, $X$ and $\Theta$ will denote vector spaces. Additionally, we assume that $X$ is a Hilbert space with inner product $\langle\cdot,\cdot\rangle$.

Consider the following equality-constrained optimization problem
\begin{subequations}\label{eq:equality-constraint-problem}
\begin{align}
    \min_{\theta \in \Theta} &\ J(\bx; \theta), \label{eq:equality-constrained-objective} \\
        \text{s.t. } & F(\bx; \theta) = 0, \label{eq:equality-constraint}
\end{align}
\end{subequations}
where the optimization parameter $\theta \in \Theta$, the state vector $\bx \in X$, the objective function $J:X\times\Theta \rightarrow \mathbb{R}$, and the equality constraint $F: X \times \Theta \rightarrow X$. 

Analogous to the unconstrained case discussed in \Cref{sec:unconstrained}, we are interested in a nonlinear splitting of the gradient for this problem, We denote the gradient of this equality-constrained problem by $\nabla_\theta J(\bx; \theta)$ where $\bx$ satisfies $F(\bx;\theta)=0$. More precisely, suppose that there exists a parameter-to-state (also known as control-to-state) mapping $\theta \mapsto \bx[\theta]$ such that $\bx[\theta]$ is a solution to $F(\bx[\theta]; \theta) = 0$. Then, the gradient is
$$ \nabla_\theta J(\bx; \theta) = \frac{d}{d\theta} J(\bx[\theta]; \theta). $$
A standard way to compute this gradient, without needing such a parameter-to-state mapping, is through the adjoint method. Let us briefly recall the standard adjoint method \cite{RoNaSa2018, GiSu2002, KuCr2016, LiPe2004, GiPi2000, const_opt} for the equality-constrained optimization problem \eqref{eq:equality-constraint-problem}. The adjoint approach computes the gradient by solving the following system
\begin{subequations}\label{eq:true-adjoint-system}
    \begin{align}
        F(\bx; \theta ) &= 0, \\
        D_1 F(\bx;\theta)^* \boldsymbol{\lambda} &= -\frac{\partial J(\bx;\theta)}{\partial \bx},\label{eq:standard-adjoint-equation} \\
        \nabla_\theta J(\bx; \theta) &= \frac{\partial J(\bx;\theta)}{\partial\theta} + \left\langle \boldsymbol{\lambda}, \frac{\partial F(\bx; \theta)}{\partial \theta} \right\rangle, \label{eq:standard-adjoint-gradient}
    \end{align}
\end{subequations}
    where $\partial/\partial\theta$ denotes the derivative with respect to explicit dependence on $\theta$, and the adjoint variable $\boldsymbol{\lambda} \in X$\footnote{The adjoint variable is more naturally interpreted as an element of the dual space $X^*$; we use the inner product on $X$ throughout to identify the dual space with $X$. For a more detailed discussion of duality in the context of adjoint equations, see \cite{TrSoBlOl2025, TrSoLe2024}.} satisfies the adjoint equation \eqref{eq:standard-adjoint-equation}, $(\cdot)^*$ denotes the adjoint with respect to the inner product.

We will now introduce the notions of nonlinear splittings of the equality constraint and of the gradient in the equality-constrained setting, and provide a systematic construction of the nonlinear splitting using a modified adjoint approach.

\begin{definition}[Nonlinear splitting of an equality constraint]\label{def:ns-grad}
A \textbf{nonlinear splitting of the equality constraint} is a mapping
\begin{align}\label{eq:nonlinear-splitting-equality-constraint}
    \mF : X \times X \times \Theta &\rightarrow X, \\
    \mF(\bx,\by; \theta) &\in X, \nonumber
\end{align}
such that $\mF(\bx,\bx; \theta) = F(\bx; \theta)$ for all $\bx \in X$ and $\theta \in \Theta$. In particular, this implies
\begin{align}\label{eq:nonlinear-splitting-equality-constraint-condition}
    \mF(\bx,\bx; \theta) = 0 \textup{ if and only if } F(\bx; \theta) = 0.
\end{align}
For each $\theta \in \Theta$, we define $K^\mF(\theta)$ as the zero level set of $\mF(\cdot,\cdot;\theta)$, i.e.,
$$ K^\mF(\theta)  \coloneqq  \{(\bx,\by) \in X \times X : \mF(\bx,\by,\theta) = 0 \}. $$
\end{definition}
We now define the notion of a nonlinear splitting of the gradient for the equality-constrained optimization problem \eqref{eq:equality-constraint-problem}. 

\begin{definition}[Nonlinear splitting of the gradient for problem \eqref{eq:equality-constraint-problem}]
    Given a nonlinear splitting of the equality constraint $\mF$, a corresponding \textbf{nonlinear splitting of the gradient} for the problem \eqref{eq:equality-constraint-problem} is given by a mapping
\begin{align}\label{eq:nonlinear-splitting-gradient-def}
    \wtnabla_\theta J(\theta) : K^{\mF}(\theta) &\rightarrow \Theta, \\
    \wtnabla_\theta J(\theta)[\bx,\by] &\in \Theta \nonumber
\end{align}
such that, if $\bx=\by$, then
\begin{equation}\label{eq:equality-gradient-diagonal}
    \wtnabla_\theta J(\theta)[\bx,\bx] = \nabla_\theta J(\bx; \theta).
\end{equation}
Again, $\nabla_\theta J(\bx; \theta)$ denotes the gradient of problem \eqref{eq:equality-constraint-problem} where $\bx$ satisfies $F(\bx;\theta)=0$.
\end{definition}

    The intuition behind the above definitions is as follows. For standard gradient-based equality-constrained optimization, for a given $\theta$, the derivative corresponding to the problem \eqref{eq:equality-constraint-problem} is evaluated at a solution to the equality constraint $\bx \in X$ to $F(\bx; \theta) = 0$.
    By splitting the equality constraint, we have enlarged the space on which we can evaluate the ``gradient''. Hence, given a nonlinear splitting of the gradient, we can consider descent algorithms which update the state vector $\bx$ and the parameters $\theta$, without requiring that the state vector solves the equality constraint at each iteration. Note that when $\bx=\by$, the equality constraint $F(\bx;\theta)=0$ is satisfied if and only if $(\bx,\bx) \in K^\mF(\theta)$. We require \eqref{eq:equality-gradient-diagonal} to further ensure consistency of the split gradient, reducing to the true gradient if $\bx=\by$.

\textbf{Nonlinear splitting from a nonlinear solver.}
    It is often the case for high-dimensional constrained optimization (e.g., PDE-constrained optimization) that solving the equality constraint at each optimization iteration is a major computational bottleneck. A broad class of numerical methods for solving such equality constraints can be interpreted as nonlinear solvers, as discussed in \cite{nonlinsolver2015}, which consist of applying an iterative fixed-point method to find the solution of the equality constraint. To mitigate this computational bottleneck, we discuss how a single iteration of such a nonlinear solver gives rise to a nonlinear splitting and hence, yields an efficient modified adjoint method to compute the gradient for the constrained objective. 
    
    Consider the solution of the equality constraint $F(\bx;\theta)=0$ via a nonlinear solver. The nonlinear solver framework of \cite{nonlinsolver2015} is a general framework which includes many common iterative methods for solving nonlinear equations, including, for example, fixed-point iterations, line searches, nonlinear multigrid, and Newton--Krylov methods. For a given $\theta$, we express one iteration of the nonlinear solver as $\bx^{k+1} = G(\bx^k; \theta)$. We assume that a solution of the original equality constraint corresponds to a fixed-point of the nonlinear solver $G$, which covers most of the examples of nonlinear solvers discussed in \cite{nonlinsolver2015}. We emphasize that $G$ is abstract and very general as an operator -- it can be include complex linear and nonlinear relations and solves within its action, and simply represents the output of a single iteration of the nonlinear fixed-point solver $G$. Note that this also includes the limiting case of $G$ being an exact solve of the equality constraint, wherein the nonlinearly split method reduces to a classical adjoint. Thus, we can redefine the equality constraint to be $F(x; \theta) := G(x; \theta) - x$. Then, define a nonlinear splitting of the equality constraint $\mF$ of $F$ by
    \begin{equation}\label{eq:nonlinear-solver-splitting}
        \mF(\bx,\by; \theta) = G(\bx;\theta) - \by.
    \end{equation}
    We modify the adjoint system \eqref{eq:true-adjoint-system} by using a single iteration of the nonlinear solver corresponding to $G$ as the state variable data for the adjoint equation. It is straightforward to verify that this defines a nonlinear splitting of the equality constraint. In particular, this gives the following nonlinear splitting adjoint gradient (NS-Adj-Gradient):
    \begin{subequations}
    \begin{align}
        \bx^{k+1} &= G(\bx^k; \theta^k), \\
        D_1 F(\bx^{k+1};\theta^k)^* \boldsymbol{\lambda}^{k+1} &= -\frac{\partial J(\bx^{k+1};\theta^k)}{\partial \bx}, \\
        \wtnabla_\theta J(\theta^k)[\bx^k,\bx^{k+1}] &=  \frac{\partial J(\bx^{k+1};\theta^k)}{\partial\theta} + \left\langle \boldsymbol{\lambda}^{k+1}, \frac{\partial F(\bx^{k+1}; \theta^k)}{\partial \theta} \right\rangle, \\
        \theta^{k+1} &= \theta^k - \gamma_k \wtnabla_\theta J(\theta^k)[\bx^k,\bx^{k+1}],
    \end{align}
    \end{subequations}
    where at each iteration, the equations are computed from top to bottom. Note that the nonlinear splitting of the gradient is subtle, with the term $\partial F/\partial \theta$ depending on both $\bx^k$ and $\bx^{k+1}$ because $\bx^{k+1}$ does not satisfy the equality constraint exactly. 
    This algorithm is summarized in \Cref{alg:ns-fp-adj-grad}.

\renewcommand{\algorithmicrequire}{\textbf{Input:}}
\renewcommand{\algorithmicensure}{\textbf{Output:}}
\begin{algorithm}[H]
\caption{NS-Adj-Gradient($\bx$, $\theta$, $\gamma$) }\label{alg:ns-fp-adj-grad}
\begin{algorithmic}
\Require state $\bx \in X$, iterate $\theta \in \Theta$, stepsize $\gamma > 0$
\State $\bx^+ = G(x; \theta)$
\State $\boldsymbol{\lambda} = $ Solve: $D_1 F(\bx^+;\theta)^* \boldsymbol{\lambda} = -\frac{\partial J(\bx^+;\theta)}{\partial \bx}$
\State $ \wtnabla_\theta J(\theta)[\bx,\bx^+] =  \frac{\partial J(\bx^+;\theta)}{\partial\theta} + \langle \boldsymbol{\lambda}, \frac{\partial F(\bx^+; \theta)}{\partial \theta} \rangle$
\State $\theta^+ = \theta - \gamma \wtnabla_\theta J(\theta)[\bx,\bx^+]$
\Ensure $[\bx^+, \theta^+, \gamma\wtnabla_\theta J(\theta)[\bx,\bx^+]]$
\end{algorithmic}
\end{algorithm}

\begin{remark}
    \Cref{alg:ns-fp-adj-grad} returns the updated state, updated iterate and the nonlinearly split gradient; we utilize all outputs when they are wrapped into accelerated methods, discussed in \Cref{sec:acceleration}.
\end{remark}

More generally, one can consider any nonlinear splitting (not necessarily one defined from a nonlinear solver). Given any nonlinear splitting of the equality constraint, $\mF$, we define a nonlinear splitting of the gradient by
    \begin{subequations}
    \begin{align}
        \mF(\bx,\by;\theta)&=0, \\
        \wtnabla_\theta J(\theta)[\bx,\by] &=  \frac{\partial J(\by;\theta)}{\partial\theta} + \left\langle \boldsymbol{\lambda}, \frac{\partial \mF(\bx,\by; \theta)}{\partial \theta} \right\rangle, \label{eq:nonlinear-splitting-gradient-adjoint-ver} \\
        \big(D_1 \mF(\bx,\by;\theta) + D_2\mF(\bx,\by;\theta)\big)^* \boldsymbol{\lambda} &= -\frac{\partial J(\by;\theta)}{\partial \by}. \label{eq:nonlinearly-split-adjoint-equation}
    \end{align}
    \end{subequations}
    It is straightforward to verify that this construction satisfies the definition of a nonlinear splitting of the gradient (\Cref{def:ns-grad}). To see this, if $\by=\bx$, we have $0 = \mF(\bx,\bx; \theta) = F(\bx; \theta)$, and the modified adjoint equation for the gradient \eqref{eq:nonlinearly-split-adjoint-equation} reduces to the standard adjoint equation \eqref{eq:standard-adjoint-equation}, since $\mF(\bx,\bx;\theta) = F(\bx; \theta)$ implies
    $$ D_1 \mF(\bx,\bx;\theta) + D_2\mF(\bx,\bx;\theta) = D_1F(\bx;\theta).$$
Thus, given a nonlinear splitting of the equality constraint, we have shown how to obtain a nonlinear splitting of the gradient using a modified adjoint method. 

\subsubsection{Acceleration}\label{sec:acceleration}
We consider accelerated gradient descent methods to wrap around the nonlinearly split adjoint gradient \Cref{alg:ns-fp-adj-grad}, using Nesterov accelerated gradient \cite{Nes83, Nes04} and Anderson acceleration (AA) \cite{AA1, AA2, Saad_2025}.

For comparison, we first state the adjoint gradient descent algorithm, modified to use the nonlinear split gradient of \Cref{alg:ns-fp-adj-grad}. We now must have two stopping criteria. The first is that the objective is sufficiently small, analogous to standard GD. The second stopping criterion enforces that the equality constraint is satisfied to sufficient tolerance by checking the residual of the equality constraint is sufficiently small. This is done by noting at the $k^{th}$ iterate of NS-GD, we have
$$ \|\delta \bx^k\| = \|\bx^{k+1} - \bx^k\| = \|F(\bx^k,\theta^k)\|. $$

\begin{algorithm}[H]
\caption{NS-Adj-GD($\bx$, $\theta$, $\eta$, $\epsilon$) }\label{alg:ns-fp-adj-gd}
\begin{algorithmic}
\Require initial state $\bx \in X$, initial iterate $\theta \in \Theta$, acceptable objective value $\eta \in \mathbb{R}$, tolerance $\epsilon > 0$
\While{$J(\bx; \theta) > \eta$ or $\|\delta \bx\| > \epsilon$}
    \State Select $\gamma > 0$
    \State $\delta \bx = \bx$
    \State $[\bx,\theta, \cdot] = $ NS-Adj-Gradient($\bx$, $\theta$, $\gamma$)
    \State $\delta \bx \mathrel{-}= \bx$
\EndWhile
\Ensure $[\bx,\theta]$
\end{algorithmic}
\end{algorithm}

The NAG variant of the nonlinear splitting adjoint gradient descent is given in \Cref{alg:ns-fp-adj-nag}.

\begin{algorithm}[H]
\caption{NS-Adj-NAG($\bx$, $\theta$, $\eta$, $\epsilon$, $\mu$) }\label{alg:ns-fp-adj-nag}
\begin{algorithmic}
\Require initial state $\bx \in X$, initial iterate $\theta \in \Theta$, acceptable objective value $\eta \in \mathbb{R}$ 
\Require tolerance $\epsilon > 0$, momentum factor $\mu > 0$
\State $\bv_{\textup{prev}} = 0$
\While{$J(\bx; \theta) > \eta$ or $\|\delta \bx\| > \epsilon$}
    \State $\theta^+ = \theta - \mu \bv_{\textup{prev}} $
    \State Select $\gamma > 0$
    \State $\delta \bx = \bx$
    \State $[\bx,\cdot, \bv] = $ NS-Adj-Gradient($\bx$, $\theta^+$, $\gamma$)
    \State $\delta \bx \mathrel{-}= \bx$
    \State $\bv_{\textup{prev}} = \mu \bv_{\textup{prev}} + \bv$
    \State $\theta = \theta - \bv_{\textup{prev}}$
\EndWhile
\Ensure $[\bx,\theta]$
\end{algorithmic}
\end{algorithm}

The AA variant is given by viewing the gradient descent step
$$ \theta^{k+1} = \theta^k - \gamma_k \wtnabla_\theta J(\theta^k)[\bx^k,\bx^{k+1}],$$
as a fixed-point iteration. A fixed-point then corresponds to a vanishing (nonlinear splitting) gradient. The AA variant is given in \Cref{alg:ns-fp-adj-AA}.

\begin{algorithm}[H]
\caption{NS-Adj-AA($\bx^0$, $\theta^0$, $\eta$, $\epsilon$, $m$) }\label{alg:ns-fp-adj-AA}
\begin{algorithmic}
\Require initial state $\bx^0 \in X$, initial iterate $\theta^0 \in \Theta$, acceptable objective value $\eta \in \mathbb{R}$
\Require tolerance $\epsilon > 0$, integer $m>0$
\State Select $\gamma > 0$
\State $[\bx^1,\theta^1,-g^0] = $ NS-Adj-Gradient($\bx^0$, $\theta^0$, $\gamma$) 
\State $[\cdot,\cdot,-g^1] = $ NS-Adj-Gradient($\bx^1$, $\theta^1$, $\gamma$)
\State $G = g^1-g^0$
\State $S = \theta^1 - \theta^0$
\For{$k=1,2,\dots$ until $J(\bx^{k}; \theta^{k}) < \eta$ and $\|\delta \bx\| < \epsilon$}
    \State $m_k = \min\{m,k\}$
    \State $\delta \bx = \bx^k$
    \State $\xi = $ QRSolve($G\xi = g^k$) \Comment{See \Cref{rmk:AA-QR}}
    \State $\theta^{k+1} = \theta^k + g^k + (S+G)\xi$
    \State $[x^{k+1},\cdot, -g^{k+1}] = $ NS-Adj-Gradient($\bx^k$, $\theta^{k+1}$, $\gamma$)
    \State $G = [G, g^{k+1}-g^k]$
    \State $S = [S, \theta^{k+1}-\theta^k]$
    \If{cols$(G)>m_k$}
        \State Remove first column of $G$ and of $S$
    \EndIf
    \State $\delta \bx \mathrel{-}= \bx^{k+1}$
\EndFor
\Ensure $[\bx,\theta]$
\end{algorithmic}
\end{algorithm}

\begin{remark}\label{rmk:AA-QR}
    In Anderson acceleration, it is common to use a (economy) QR decomposition to solve $G\xi = g^k$, interpreted in the least squares sense. Note that the QR decomposition does not need to be recomputed at each iteration, but rather updated when a column of $G$ is added or removed. For more details, see \cite{AA1, AA2, Saad_2025}.
\end{remark}

In the third numerical example, \Cref{sec:neutron-transport}, we will compare these methods to their standard counterparts, for a source optimization problem where the equality constraint is given by a nonlinear neutron transport equation.

\section{Numerical examples}\label{sec:numerical-example}
We consider several numerical examples that demonstrate the nonlinear splitting framework. In \Cref{sec:uncons-examples}, we demonstrate nonlinear splitting methods on examples in unconstrained optimization; particularly, in \Cref{sec:benchmarks}, we consider several benchmark nonconvex functions and in \Cref{sec:NLS}, we consider a nonlinear least squares problem arising in pure state quantum tomography. Subsequently, in \Cref{sec:cons-examples}, we demonstrate the nonlinear splitting methods on examples in equality-constrained optimization; particularly, we consider the stationary neutron transport equation as the equality constraint in \Cref{sec:neutron-transport} and an initial value problem as the equality constraint in \Cref{sec:time-int-example}.

For convenience, we summarize the methods introduced thus far in \Cref{table:list-methods}.

\begin{table}[H]
\centering
\begin{tabular}{|l|l|l|l|} 
 \hline
 Name & Optimization type & Method & Standard counterpart \\ [0.5ex] 
 \hline
 NS-GD & Unconstrained & \eqref{eq:GD-partitioned} & GD \eqref{eq:GD-standard} \\ 
 NS-NAG & " & \eqref{eq:NAG-partitioned} & NAG \eqref{eq:NAG-standard} \\
 NS-GD-Newton$(1)$ & " & \eqref{eq:GD-splitting-Newton-1} & GD \\
 NS-NAG-Newton$(1)$ & " & \eqref{eq:NS-NAG-Newton-1} & NAG \\
 NS-AA-Newton$(1)$ & " & \eqref{eq:NS-AA-Newton-1} & AA applied to GD \\ 
  NS-Adj-GD & Equality-constrained & \Cref{alg:ns-fp-adj-gd} & Adj-GD (see \Cref{sec:neutron-transport}) \\ 
  NS-Adj-NAG & " & \Cref{alg:ns-fp-adj-nag} & Adj-NAG (see \Cref{sec:neutron-transport})\\ 
  NS-Adj-AA & " & \Cref{alg:ns-fp-adj-AA} & Adj-AA (see \Cref{sec:neutron-transport}) \\ 
 \hline
\end{tabular}
\caption{List of nonlinear splitting methods for unconstrained and equality-constrained optimization}
\label{table:list-methods}
\end{table}

\subsection{Unconstrained examples}\label{sec:uncons-examples}

\subsubsection{Benchmark nonconvex optimization problems}\label{sec:benchmarks}
\begin{figure}[!ht]
\subfloat[Rastrigin]{\includegraphics[trim={0.5cm 0.2cm 0.5cm 1.35cm},clip,width = 2.2in]{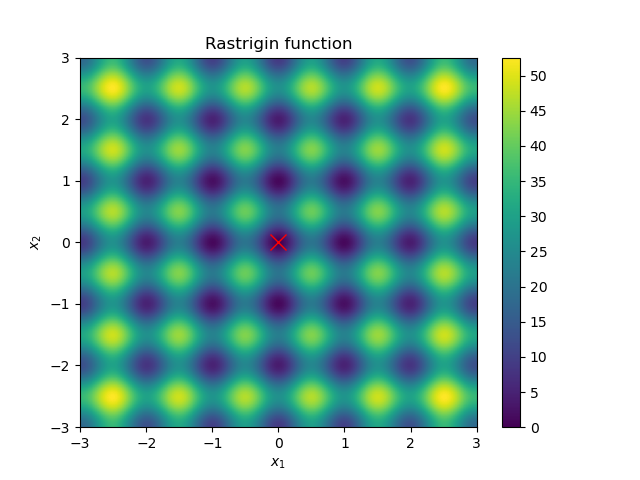}}
\subfloat[Rosenbrock]{\includegraphics[trim={0.5cm 0.2cm 0.5cm 1.35cm},clip,width = 2.2in]{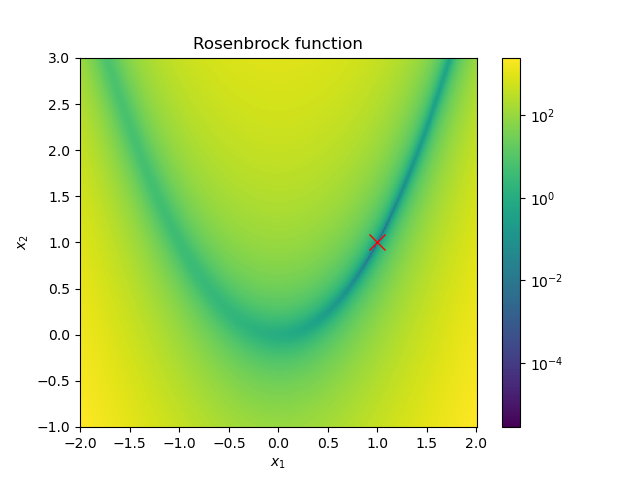}}
\subfloat[Beale]{\includegraphics[trim={0.5cm 0.2cm 0.5cm 1.35cm},clip,width = 2.2in]{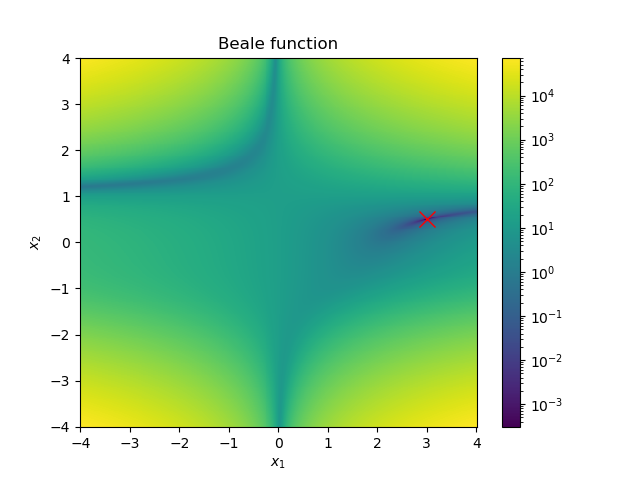}}
\caption{Plots of the three benchmark functions. Red `x' marks the global minimum.}\label{fig:func-plots}
\end{figure}

As examples of challenging unconstrained non-convex optimization using nonlinear splitting, we consider several two-dimensional non-convex benchmark problems, namely, the Rastrigin, Rosenbrock, and Beale functions, which all have a global minimum value of zero, and complex nonconvex loss landscapes. The functions and their global minima are shown in \Cref{fig:func-plots}.
For each of these functions, we define a linearly-implicit nonlinear splitting (\Cref{ex:li-gradient}), and we compare NS-GD and NS-NAG based on these splittings to GD, NAG, discrete gradient NS-GD, discrete gradient NS-NAG, and to the built-in MATLAB \cite{MATLAB} optimizers \verb|fmin| (FM), \verb|Pattern search| (PS), \verb|Genetic algorithm| (GA), \verb|Particle swarm| (PartS), \verb|Surrogate| (Surr), \verb|Global search| (GS).

The Rastrigin function, its gradient, and a linearly-implicit nonlinear splitting of the gradient are given by
\begin{align*}
    C(\bx) &= 20 + (x_1^2-10\cos(2\pi x_1)) + (x_2^2 - 10\cos(2\pi x_2)), \\
    \nabla C(\bx) &= \begin{bmatrix} 2x_1 + 20\pi \sin(2\pi x_1) \\ 2x_2 + 20\pi \sin(2\pi x_2)\end{bmatrix}, \\
    \wtnabla C(\bx,\by) &= \begin{bmatrix} 2y_1 + 20\pi \sin(2\pi x_1) y_1/x_1 \\ 2y_2 + 20\pi \sin(2\pi x_2) y_2/x_2\end{bmatrix} .
\end{align*}
The Rosenbrock function, its gradient, and a linearly-implicit nonlinear splitting of the gradient are given by
\begin{align*}
    C(\bx) &= (1-x_1)^2 + 100(x_2-x_1^2)^2, \\
    \nabla C(\bx) &= \begin{bmatrix} 2(200 x_1^3 - 200 x_1x_2 +x_1-1) \\ 200(x_2-x_1^2) \end{bmatrix}, \\
    \wtnabla C(\bx,\by) &= \begin{bmatrix} 2(200 y_1x_1^2 - 200 x_1x_2 +x_1-1) \\ 200(y_2-x_1y_1) \end{bmatrix} .
\end{align*}
The Beale function, its gradient, and a linearly-implicit nonlinear splitting of the gradient are given by
\begin{align*}
    C(\bx) &= (1.5-x_1 +x_1x_2)^2 + (2.25-x_1+x_1x_2^2)^2 + (2.625-x_1+x_1x_2^3)^2, \\
    \nabla C(\bx) &= \begin{bmatrix} 2x_1(x_2^6 + x_2^4 - 2x_2^3 - x_2^2 + 5.25 x_2^3) + 5.25x_2^3 + 4.5x_2^2 + 3x_2 - 12.75 \\  6x_1( x_1(x_2^5 + \frac{2}{3}x_2^3 - x_2^2 - \frac{1}{3}x_2 - \frac{1}{3}) + 2.625x_2^2 + 1.5 x_2 + 0.5) \end{bmatrix}, \\
    \wtnabla C(\bx,\by) &= \begin{bmatrix} 2y_1(x_2^6 + x_2^4 - 2x_2^3 - x_2^2 + 5.25 x_2^3) + 5.25x_2^2y_2 + 4.5x_2y_2 + 3y_2 - 12.75 \\  6y_1( x_1(x_2^5 + \frac{2}{3}x_2^3 - x_2^2 - \frac{1}{3}x_2 - \frac{1}{3}) + 2.625x_2^2) + 6x_1(1.5 y_2 + 0.5) \end{bmatrix}.
\end{align*}

In the following, NS-GD and NS-NAG refer to the nonlinear splitting methods with the above linearly-implicit splittings, whereas NS-GD(AVF) and NS-NAG(AVF) refer to the nonlinear splitting methods using the average vector field (\Cref{ex:disc-grad}) as the splitting.

\begin{figure}[!ht]
\subfloat[]{\includegraphics[trim={3.7cm 8.5cm 4.5cm 9.15cm},clip,width = 2.7in]{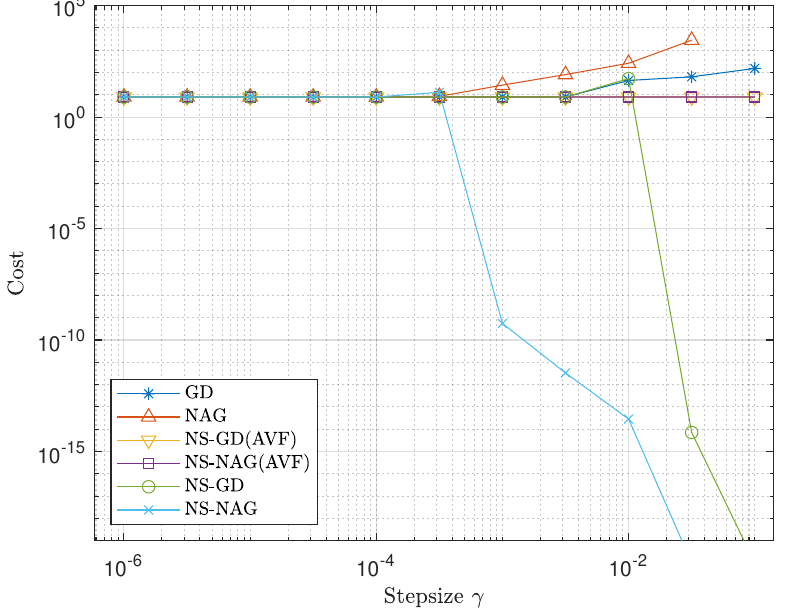}}
\subfloat[]{\includegraphics[trim={3.7cm 8.5cm 4.5cm 9.15cm},clip,width = 2.7in]{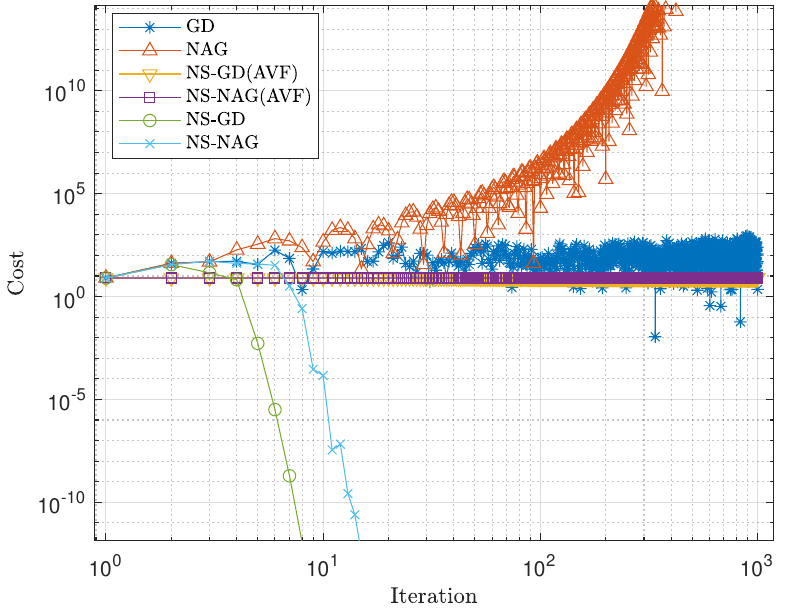}}
\caption{Rastrigin function with initial iterate $\bx^0 = (2,2)$. (A) Cost after 10,000 iterations as a function of stepsize. (B) Cost versus iteration for stepsize $\gamma = 0.1$.}\label{fig:cost-rastrigin}
\end{figure}

\begin{figure}[!ht]
\subfloat[]{\includegraphics[trim={3.7cm 8.5cm 4.5cm 9.15cm},clip,width = 2.7in]{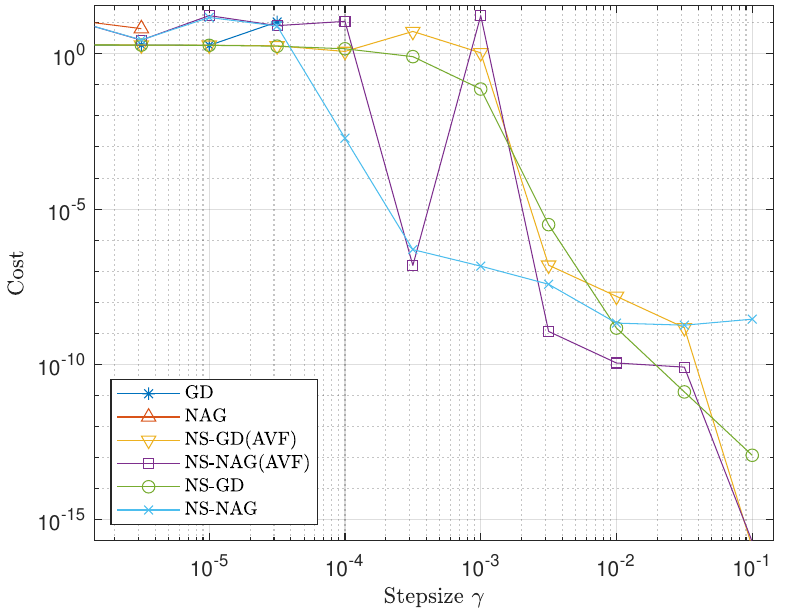}}
\subfloat[]{\includegraphics[trim={3.7cm 8.5cm 4.5cm 9.15cm},clip,width = 2.7in]{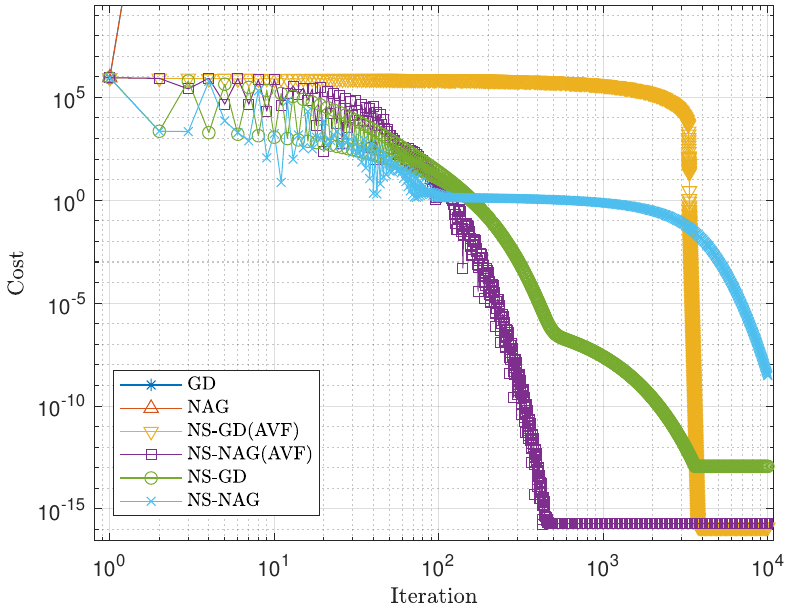}}
\caption{Rosenbrock function with initial iterate $\bx^0 = (10,5)$. (A) Cost after 10,000 iterations as a function of stepsize. (B) Cost versus iteration for stepsize $\gamma = 0.1$.}\label{fig:cost-rosenbrock}
\end{figure}

\begin{figure}[!ht]
\subfloat[]{\includegraphics[trim={3.7cm 8.5cm 4.5cm 9.15cm},clip,width = 2.7in]{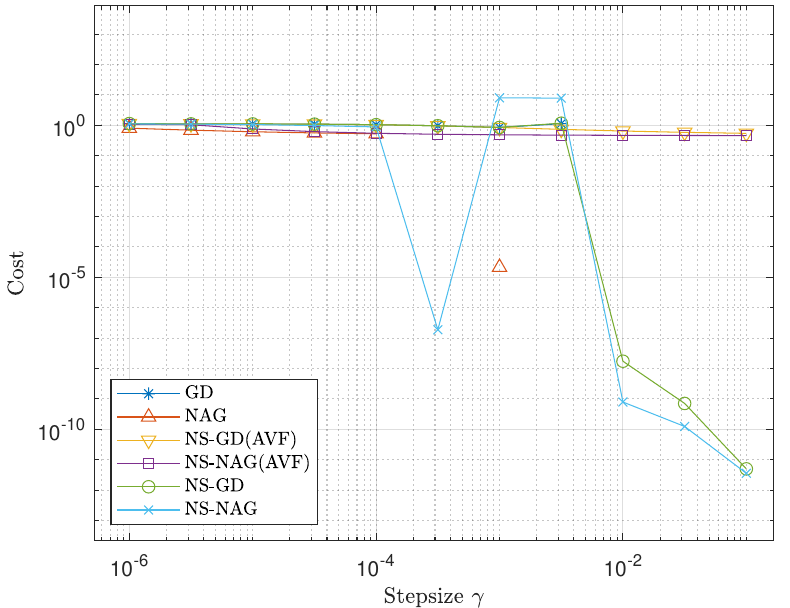}}
\subfloat[]{\includegraphics[trim={3.7cm 8.5cm 4.5cm 9.15cm},clip,width = 2.7in]{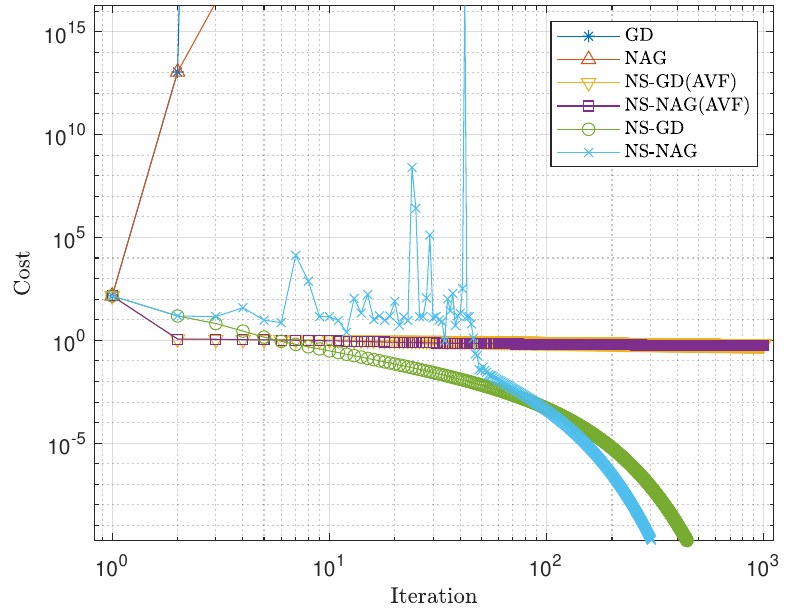}}
\caption{Beale function with initial iterate $\bx^0 = (-2,2)$. (A) Cost after 10,000 iterations as a function of stepsize. (B) Cost versus iteration for stepsize $\gamma = 0.1$.}\label{fig:cost-beale}
\end{figure}

\Cref{fig:cost-rastrigin,fig:cost-rosenbrock,fig:cost-beale} show the cost as a function of stepsize and iteration for these three different problems. For all three cases, classical GD and NAG are unable to achieve $<\mathcal{O}(1)$ loss in 10,000 iterations for any stepsizes tested across $\gamma\in[10^{-6},0.1]$ and furthermore, are unstable for larger stepsizes. In contrast, the linearly-implicit methods NS-GD and NS-NAG are stable for larger step sizes for all three problems, and are able to converge to the true minimum. The AVF variants are only successful for the Rosenbrock problem \Cref{fig:cost-rosenbrock}, although there they do converge notably faster than the linearly-implicit methods. 

\Cref{table:time-cost-matlab-comp} shows a comparison of average runtime and average cost between NS-GD, NS-NAG and built-in MATLAB optimizers. The runs are averaged over 100 random initial conditions, using the same random initial conditions for all methods. For the Rastrigin function, $\bx^0 \in [0,3] \times [0,3]$; for Rosenbrock, $\bx^0 \in [-2,2] \times [-1,3]$; for Beale, $\bx^0 \in [-4,4]\times[-4,4]$. NS-GD and NS-NAG use $\gamma = 0.1$ and are run until either 4000 iterations or the cost goes below $10^{-10}$. The most important point from \Cref{table:time-cost-matlab-comp} is that the only methods to obtain a small average cost across all three problems (i.e. find the true minimum more or less independent of initial condition) are NS-GD, NS-NAG, and the MATLAB Global Search (GS). In terms of runtime, the NS methods are $6-9\times$ faster for Rastrigin, approximately the same for Beale, and $4-5\times$ slower for Rosenbrock. That being said, note that the NS methods are implemented via naive looping in MATLAB scripts, and we expect a relatively significant performance degradation compared with the optimized methods provided in MATLAB. 

\begin{table}[!htb]
\centering
    \begin{tabular}{|l|l|l|l|l|l|l|} 
     \hline
     & \multicolumn{2}{c|}{Rastrigin} & \multicolumn{2}{c|}{Rosenbrock} & \multicolumn{2}{c|}{Beale} \\\hline 
     Method & Runtime (s) & Cost & Runtime (s) & Cost & Runtime (s) & Cost \\
     \hline
      NS-GD & 0.0064 & 4.3e-12 & 1.42 & 2.09e-9 & 0.12 & 2.27e-9 \\
      NS-NAG & 0.0089 & 3.49e-11 & 1.45 & 2.46e-5 & 0.11 & 1.07e-10 \\  
      FM & 0.0012 & 4.97 & 0.0013 & 9.85 & 0.0013 & 1.02e-9 \\
      PS & 0.018 & 3.41e-11 & 0.026 & 8.01e-10 & 0.017 & 1.03 \\
      GA & 0.010 & 1.06 & 0.012 & 1.24 & 0.0085 & 0.18 \\
      PartS & 0.0085 & 0.22 & 0.011 & 0.18 & 0.0055 & 0.0093 \\
      Surr & 1.36 & 0 & 2.61 & 0.014 & 1.59 & 0.013 \\
      GS & 0.059 & 0 & 0.31 & 0 & 0.11 & 3.93e-14 \\
     \hline
    \end{tabular} 
\caption{Comparison of runtime and cost for the three benchmark functions, averaged over 100 randomized initial states, for NS-GD, NS-NAG, and MATLAB optimizers.}
\label{table:time-cost-matlab-comp}
\end{table}


\subsubsection{Nonlinear least squares: pure state quantum tomography}\label{sec:NLS}
In this example, we consider a problem in pure state quantum tomography \cite{tomography2020}, which is to construct a state $|\Psi\rangle$ such that the expectation value of a set of operators $\{A_i\}_{i=1}^N$ on this state, $\langle \Psi| A_i |\Psi\rangle$, agree with a set of observed values $\{a_i\}_{i=1}^N$. We approach this problem as an unconstrained state reconstruction. This can be formulated as a nonlinear least squares problem of the form
\begin{align}\label{eq:nonlinear-least-squares-quantum}
    \min_{|\Psi\rangle} \sum_i ( \langle \Psi| A_i |\Psi\rangle - a_i )^2.
\end{align}
For this example, we consider $n$-qubit states
$$ |\Psi\rangle = \sum_{i_1,\dots,i_n = 0}^1 c_{i_1 \dots i_n} |i_1\rangle \otimes \cdots \otimes |i_n\rangle, $$
where $|0\rangle$ and $|1\rangle$ denote the standard single qubit ground and excited states, and $c_{i_1 \dots i_n} \in \mathbb{C}$. The dimension of this vector space is $2^n$ as a complex vector space and $2^{n+1}$ as a real vector space. A basis for the space of linear operators on this space can be constructed by considering the extended Pauli basis
\begin{equation}\label{eq:pauli-basis}
    \hat{\sigma}_0 = I = \begin{bmatrix} 1 & 0 \\ 0 & 1 \end{bmatrix}, \hat{\sigma}_1 = \begin{bmatrix} 0 & 1 \\ 1 & 0 \end{bmatrix}, \hat{\sigma}_2 = \begin{bmatrix} 0 & -i \\ i & 0 \end{bmatrix}, \hat{\sigma}_3 = \begin{bmatrix} 1 & 0 \\ 0 & -1 \end{bmatrix}.
\end{equation}
A basis for operators on $n$-qubit states is then given by $n$-fold tensor products of the extended Pauli basis, i.e.,
\begin{equation}\label{eq:pauli-basis-tensor}
     \{ \hat{\sigma}_{j_1} \otimes \cdots \otimes \hat{\sigma}_{j_n} \}_{j_1,\dots,j_n = 0}^3.
\end{equation}

\begin{figure}[!htb]
\subfloat[]{\includegraphics[trim={3.7cm 8.5cm 4.5cm 9.15cm},clip,width = 2.15in]{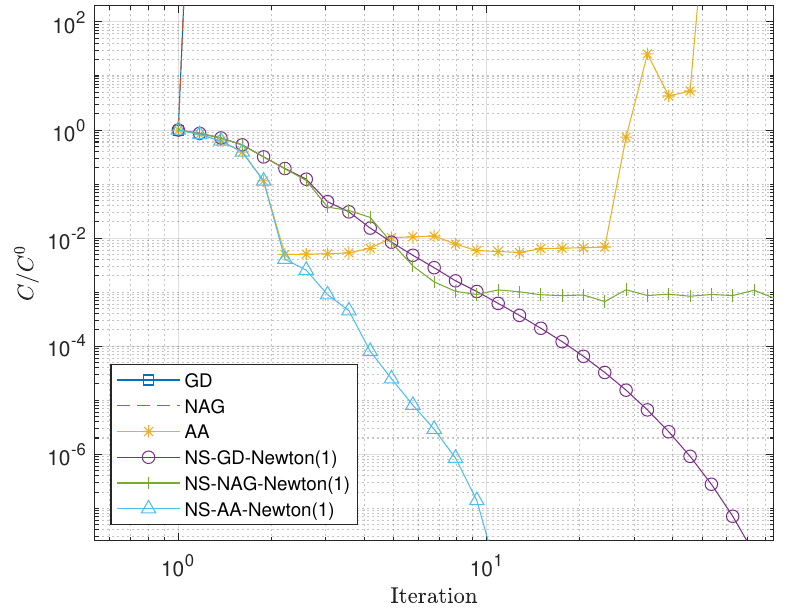}}
\subfloat[]{\includegraphics[trim={3.7cm 8.5cm 4.5cm 9.15cm},clip,width = 2.15in]{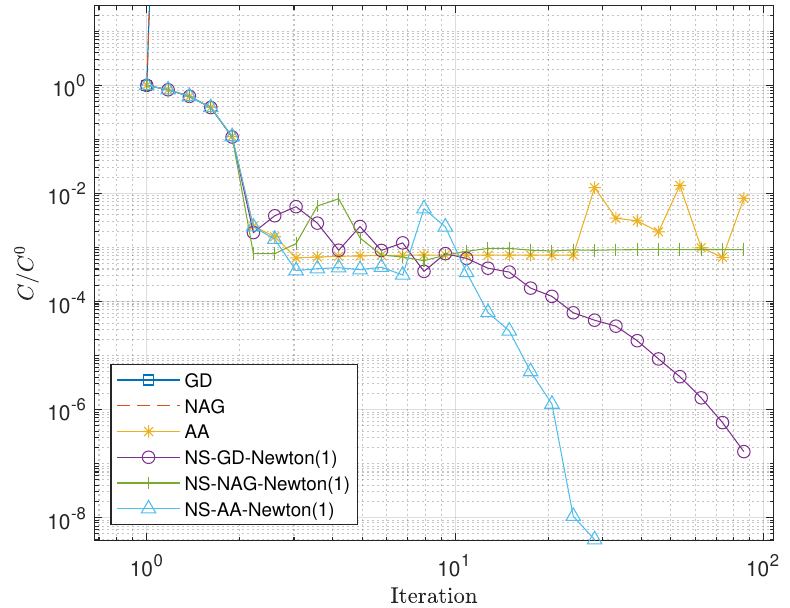}}
\subfloat[]{\includegraphics[trim={3.7cm 8.5cm 4.5cm 9.15cm},clip,width = 2.15in]{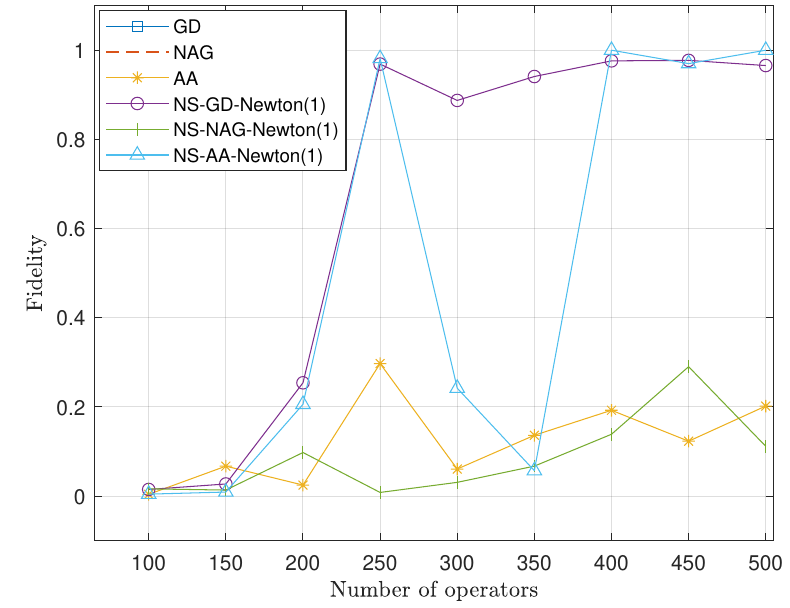}}
\caption{(A) Normalized cost versus iteration for 100 operators. (B) Normalized cost versus iteration for 500 operators. (C) Fidelity after 100 iterations versus number of operators.}\label{fig:quantum-plots}
\end{figure}

For this example, we consider $n=6$ which corresponds to a real vector space of dimension $2^{6+1} = 128$ with $4^6 = 4098$ operators defined by \eqref{eq:pauli-basis-tensor}. To generate the values $a_i$ for the problem \eqref{eq:nonlinear-least-squares-quantum}, we constructed a random quantum state $|\Psi_*\rangle$ and computed its expectation value with all 4098 operators. In \Cref{fig:quantum-plots}, we compare NS-GD-Newton(1), NS-NAG-Newton(1), and NS-AA-Newton(1), using the nonlinear least squares splitting defined in \Cref{ex:LM-splitting}, to their standard counterparts GD, NAG, and AA, all using stepsize $\gamma=10$. In \Cref{fig:quantum-plots}(A), we compare the cost versus iteration of these methods with 100 randomly chosen operators. Similarly, in \Cref{fig:quantum-plots}(B), we compare these methods for 500 randomly chosen operators. In \Cref{fig:quantum-plots}, we compare the fidelity of the optimized state after 100 iterations versus the number of operators chosen, where the fidelity of the optimized state $|\Psi\rangle$ versus the prepared state $|\Psi_*\rangle$ is defined as
\begin{align*}
    \text{Fidelity}(|\Psi\rangle,|\Psi_*\rangle) &:= \left(\text{trace} \sqrt{\sqrt{\rho(\Psi)} \rho(\Psi_*) \sqrt{\rho(\Psi)}} \right)^2, \\
    \rho(\Psi) &:= |\Psi\rangle \langle \Psi |,
\end{align*}
with a fidelity of $1$ corresponding to equivalent states. In \Cref{fig:quantum-plots}, standard GD and NAG are immediately unstable. On the other hand, all of the nonlinear splitting methods are stable. Furthermore, both NS-GD-Newton(1) and NS-AA-Newton(1) converge to the true minimum of the cost. For sufficiently many operators, the fidelity of the reconstructed state is close to 1 for both NS-GD-Newton(1) and NS-AA-Newton(1). For a lower number of operators, despite both NS-GD-Newton(1) and NS-AA-Newton(1) reaching the true minimum of the cost, the fidelity is relatively compromised due to random selection of operator subset without any systematic selection scheme.

\subsection{Examples with equality constraints}\label{sec:cons-examples}

\subsubsection{Nonlinear neutron transport equation}\label{sec:neutron-transport}
As an example of equality-constrained optimization using nonlinear splitting, we consider a source optimization problem for the nonlinear neutron transport equation. This source optimization problem arises in optical molecular imaging \cite{transportsource2} and bioluminescence tomography \cite{transportsource1}; well-posedness of the optimization problem with linear transport is investigated in \cite{transportsource1, transportsource2, transportsource3}.

The equality constraint for this optimization problem is given by the following nonlinear neutron transport equation
\begin{subequations}\label{eq:neutron-transport}
\begin{align}
    \bsO \cdot \nabla \psi(x, \bsO) + \sigma_t(x) \psi(x,\bsO) &- \frac{\sigma_s(x,\phi(x))}{4\pi} \int_{\mathbb{S}^2} \psi(x,\bsO) d\Omega = q(x,\bsO), \label{eq:neutron-transport-a}\\
    \psi(x,\bsO) &= \psi_{\textup{in}}(x,\bsO) \text{ for } x \in \partial D \text{ and } \bsO \cdot \mathbf{n} < 0,
\end{align}
\end{subequations}
where $\bsO \in \mathbb{S}^2$ is the transport angle, $d\Omega$ is the standard measure on $\mathbb{S}^2$, $D \ni x$ is the spatial domain, $q(x,\bsO)$ is the source, $\psi_{\textup{in}}$ is the inflow boundary data, $\sigma_t(x)$ is the total cross section, and $\sigma_s(x,\phi(x))$ is the nonlinear scattering cross section which depends on $x$ and the scalar flux $\phi$, defined as the zeroth angular moment of $\psi$, i.e.,
$$ \phi(x)  \coloneqq  \int_{\mathbb{S}^2} \psi(x,\bsO) d\Omega. $$
Due to the dependence of the scattering cross section on $\phi$ (and hence, $\psi$), \eqref{eq:neutron-transport} is a nonlinear equality constraint.

We consider the following equality-constrained optimization problem
\begin{align}\label{eq:source-opt-general}
    \min_{q(x,\bsO)} J(\psi) &= \frac{1}{2} \| \psi - \psi_{\text{tar}} \|_{L^2_{\bx,\bsO}}^2, \\
    &\text{with equality constraint } \eqref{eq:neutron-transport}, \nonumber
\end{align}
where the norm is induced from the inner product
$$ \left( f, g \right)_{L^2_{x,\bsO}}  \coloneqq  \int_D \int_{\mathbb{S}^2} f(x,\bsO) g(x,\bsO) d\Omega dx. $$
i.e., to find a source $q(x,\bsO)$ such that the solution to the nonlinear transport equation \eqref{eq:neutron-transport} matches a target state $\psi_{\text{tar}}$. Additional details on the adjoint equation of the neutron transport equation for computing the gradient are provided in \Cref{sec:neutron-transport-adjoint} (see also \cite{Prinja2010}).

For this example, we take $D = [0, 0.5] \times [0,1]$ using a uniform $24$-by-$48$ grid, $\sigma_t(x) = 2$ for all $x$, $\sigma_s(x,\phi(x)) = 1 + 0.02 \phi(x)^2$, and inflow boundary data
\begin{align*}
    \psi_{\text{in}}(x,\bsO) &= \frac{\pi}{4} \Big( \sin(\pi x_1) \sin(\pi x_2) \\ 
    & \qquad \qquad + (\Omega_1 + \Omega_2)\sin(2 \pi x_1) \sin(2 \pi x_2) + \frac{1}{4} (\Omega_1^2 + \Omega_2^2)\sin(3\pi x_1) \sin(3 \pi y_1) + 2 \Big),
\end{align*}
where $x = (x_1,x_2)$ denote the spatial coordinates and $\bsO = (\Omega_1,\Omega_2,\Omega_3)$ denote the angular coordinates using the embedding $\mathbb{S}^2 \hookrightarrow \mathbb{R}^3$. The target state $\psi_{\text{tar}}$ was generated by solving the nonlinear neutron transport equation \eqref{eq:neutron-transport} with the above boundary data, cross sections, and a source given by
\begin{align*}
    q(x, \bsO) &= \Omega_1 X_1(x,\bsO) + \Omega_2 X_2(x,\bsO) + 2\psi_{\text{in}}(x,\bsO) - \frac{1}{4\pi} X_3(x), \\
    X_1(x,\bsO) & \coloneqq  \frac{1}{4}\cos(\pi x_1) \sin(\pi x_2) + \frac{1}{2}(\Omega_1 + \Omega_2) \cos(2\pi x_1) \sin(2 \pi x_2), \\
    X_2(x,\bsO) & \coloneqq  \frac{1}{4}\sin(\pi x_1) \cos(\pi x_2) + \frac{1}{2}(\Omega_1 + \Omega_2) \sin(2\pi x_1) \cos(2 \pi x_2), \\
    X_3(x) & \coloneqq  \sin(\pi x_1) \sin(\pi x_2) + \frac{1}{6} \sin(3\pi x_1) \sin(3 \pi x_2) + 2.
\end{align*}
The spatial discretization is constructed with the finite element library MFEM \cite{mfem,mfem-web} using linear $L^2$ discontinuous finite elements, with $ 24 \times 48 \times 4 = 4608$ spatial degrees of freedom. The angular dependence is discretized using the $S_N$ method (see, e.g., \cite{CaoWu2020, Thynell1998}) with the Level Symmetric $S_{20}$ angular quadrature rule which has 220 discrete directions, yielding $4608 \times 220 = 1013760$ degrees of freedom for functions of space and angle, particularly for the state vector $\psi$ and the optimization parameter (source) $q$.

The nonlinear splitting of the equality constraint is given by linearizing the scattering cross section about the previous optimization iterate, i.e.,
\begin{align}
    \bsO \cdot \nabla \psi^{k+1}(x, \bsO) + \sigma_t(x) \psi^{k+1}(x,\bsO) &- \frac{\sigma_s(x,\phi^k(x))}{4\pi} \int_{\mathbb{S}^2} \psi^{k+1}(x,\bsO) d\Omega = q^{k}(x,\bsO), \label{eq:neutron-transport-split}
\end{align}

\begin{figure}[!htb]
\subfloat[]{\includegraphics[trim={3.7cm 8.5cm 4.5cm 9.15cm},clip,width = 2.15in]{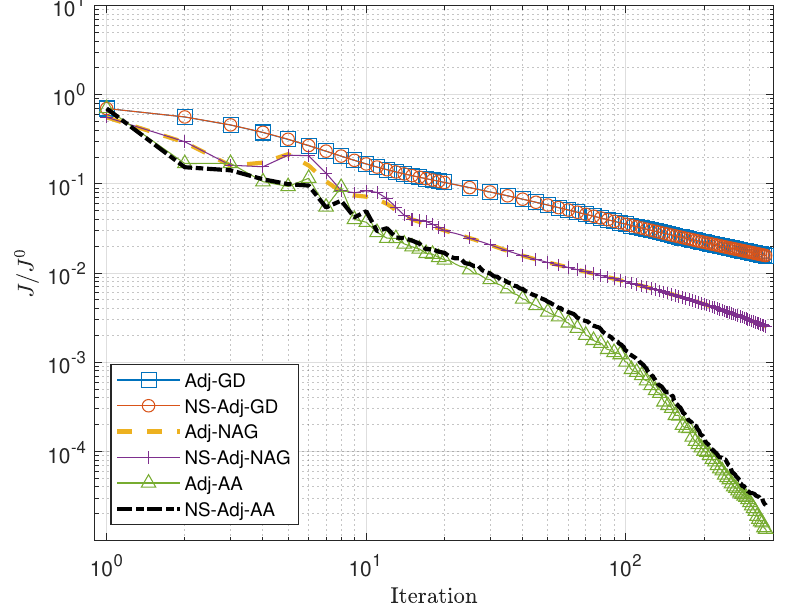}}
\subfloat[]{\includegraphics[trim={3.7cm 8.5cm 4.5cm 9.15cm},clip,width = 2.15in]{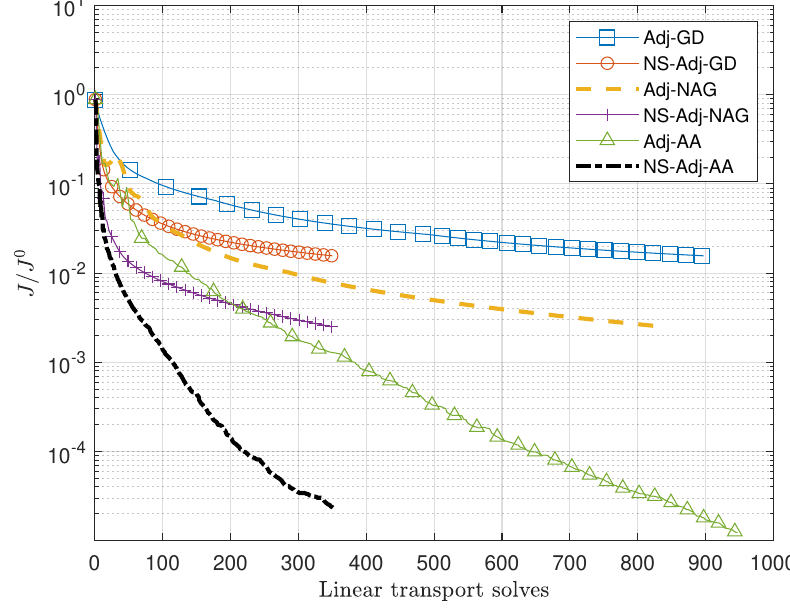}}
\subfloat[]{\includegraphics[trim={3.7cm 8.5cm 4.5cm 9.15cm},clip,width = 2.15in]{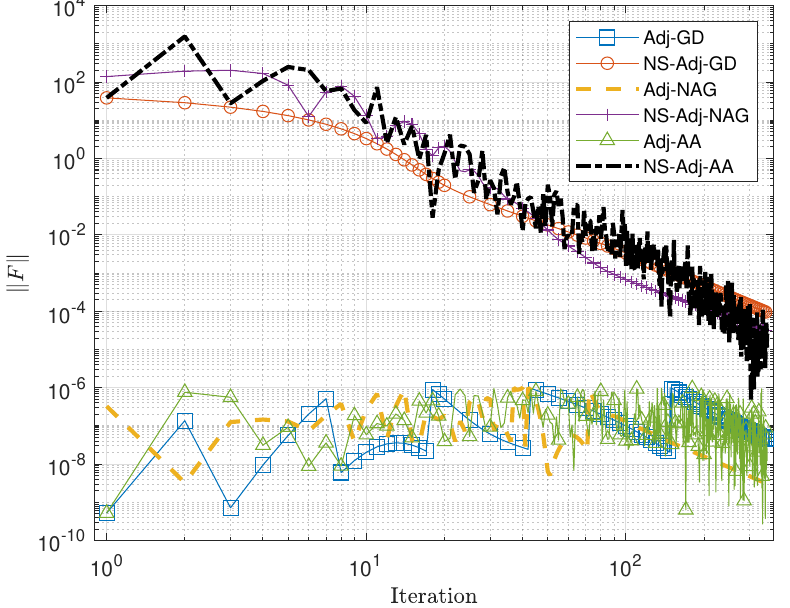}}
\caption{(A) Relative objective versus iteration. (B) Efficiency comparison: relative objective versus linear transport solves in forward pass. (C) Equality constraint residual versus iteration }\label{fig:npg-obj}
\end{figure}

We compare NS-Adj-GD (\Cref{alg:ns-fp-adj-gd}), NS-Adj-NAG (\Cref{alg:ns-fp-adj-nag}), NS-Adj-AA (\Cref{alg:ns-fp-adj-AA}) with their standard counterparts Adj-GD, Adj-NAG, and Adj-AA, where the standard methods solve the nonlinear equality constraint to a tolerance of $10^{-6}$ at each optimization iteration. A constant $\gamma = 5 \times 10^{-9}$ is used for all optimization methods. NS-Adj-NAG and Adj-NAG both use a momentum factor $\mu=0.9$. NS-Adj-AA and Adj-AA both use $m=25$. All runs are initialized with zero initial state, $\psi^0 = 0$, and zero initial source, $q^0 = 0$. \Cref{fig:npg-obj}(A) shows the relative objective $J/J^0$ versus iteration, where $J^0$ is the objective value with the initial state and source $\psi^0 = 0 = q^0$. An efficiency plot is shown in \Cref{fig:npg-obj}(B), which compares the relative objective versus number of linear transport solves in the forward pass; in particular, the nonlinear splitting variants are able to achieve the same cost as the standard methods with $2-3\times$ less linear transport solves in the forward pass. The equality constraint residual $\|F\|$ versus iteration is shown in \Cref{fig:npg-obj}(C). 

\begin{figure}[!htb]
\stackunder[1 pt]{\includegraphics[trim={16cm 2cm 21cm 2cm},clip,width = 1.655in]{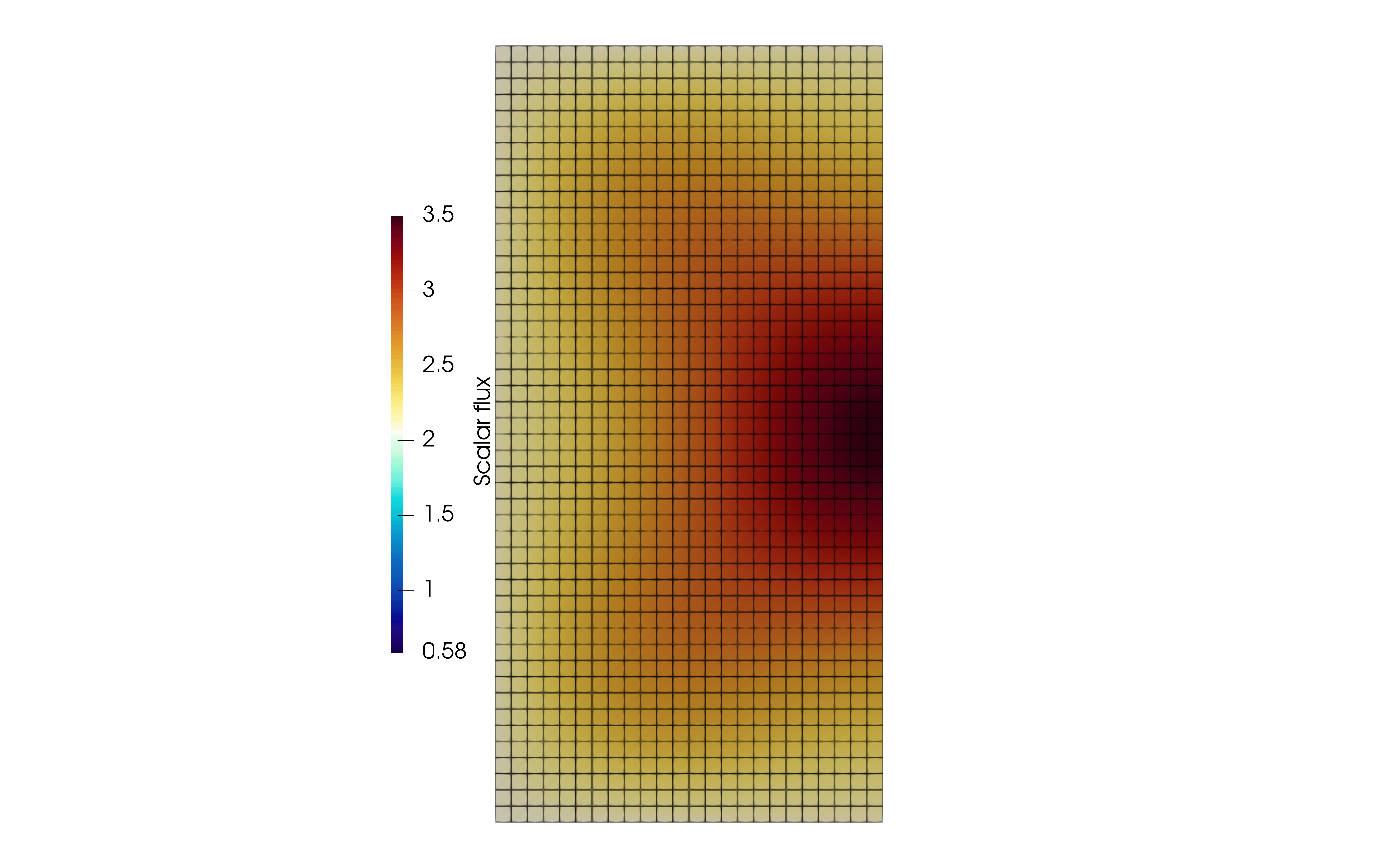}}{\qquad $\phi_{\text{tar}}$}
\stackunder[1 pt]{\includegraphics[trim={21cm 2cm 21cm 2cm},clip,width = 1.35in]{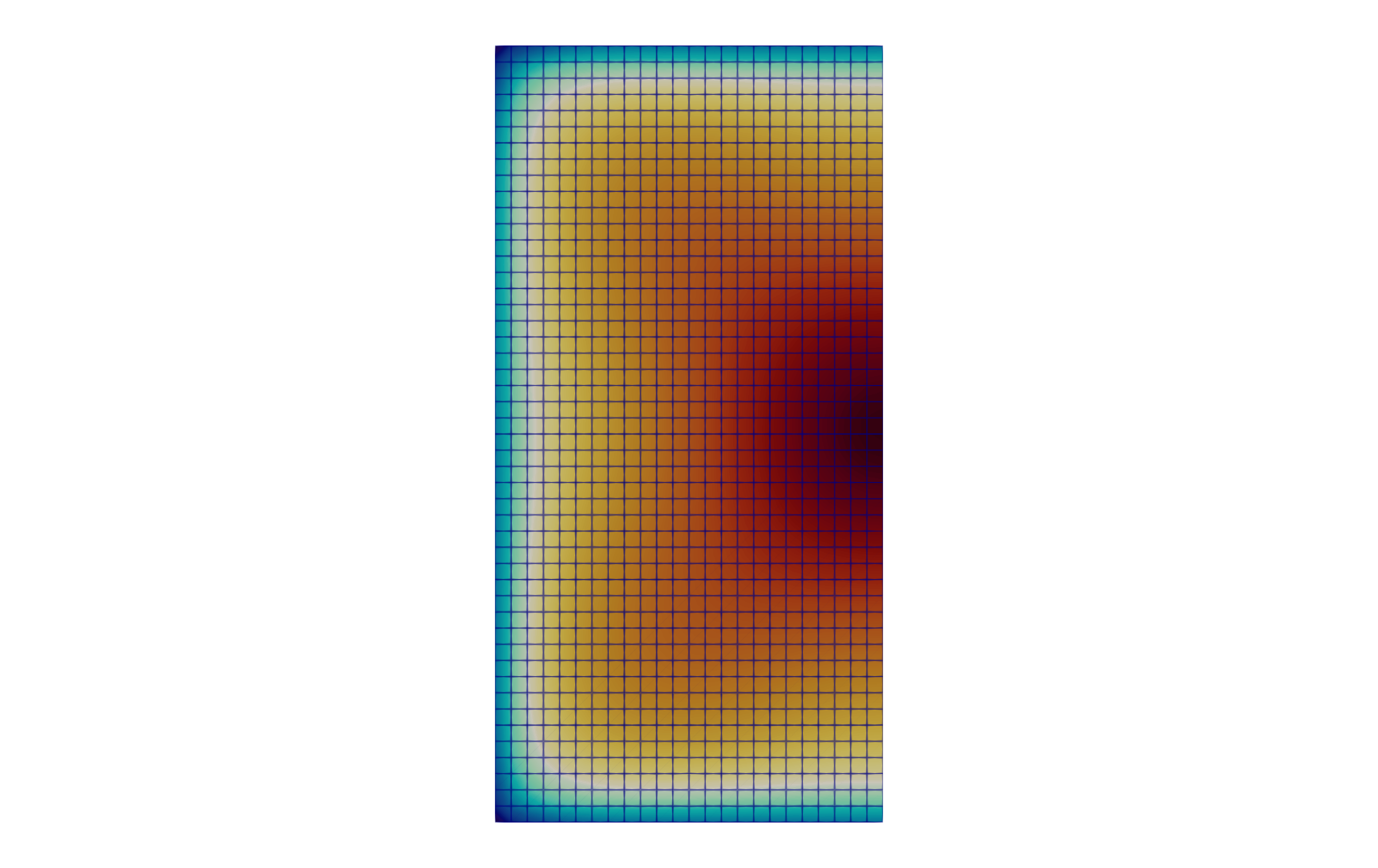}}{$\phi_{\text{NS-Adj-GD}}$}
\stackunder[1 pt]{\includegraphics[trim={21cm 2cm 21cm 2cm},clip,width = 1.35in]{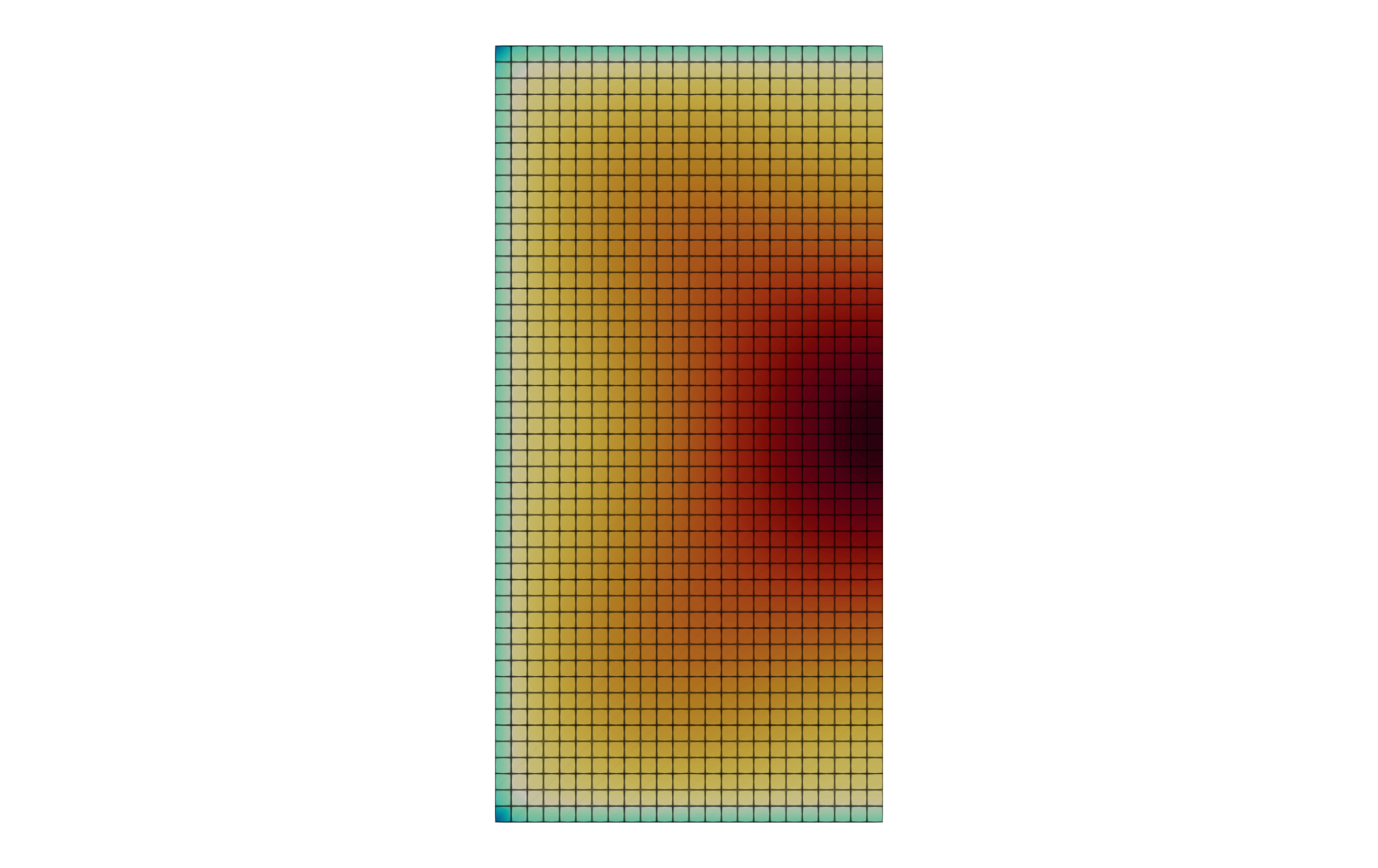}}{$\phi_{\text{NS-Adj-NAG}}$}
\stackunder[1 pt]{\includegraphics[trim={21cm 2cm 21cm 2cm},clip,width = 1.35in]{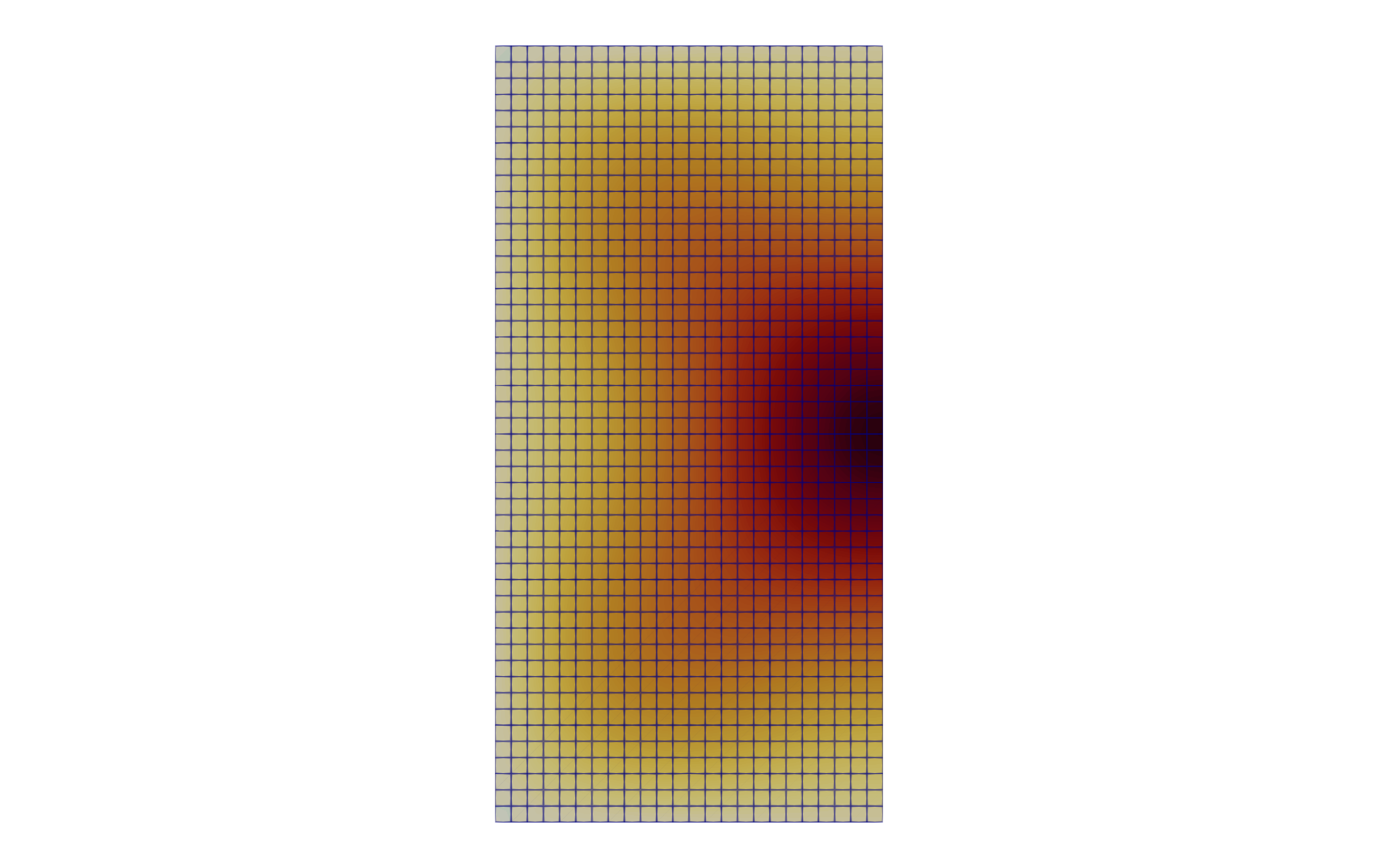}}{$\phi_{\text{NS-Adj-AA}}$}
\caption{Comparison of the scalar flux of the target state ($\phi_{\text{tar}}$) versus states optimized using NS-Adj-GD ($\phi_{\text{NS-Adj-GD}}$), NS-Adj-NAG ($\phi_{\text{NS-Adj-NAG}}$), and NS-Adj-AA ($\phi_{\text{NS-Adj-AA}}$). All plots use the same scale.}\label{fig:phi-solution-comp}
\end{figure}

\begin{figure}[!htb]
\stackunder[1 pt]{\includegraphics[trim={16cm 2cm 21cm 2cm},clip,width = 1.655in]{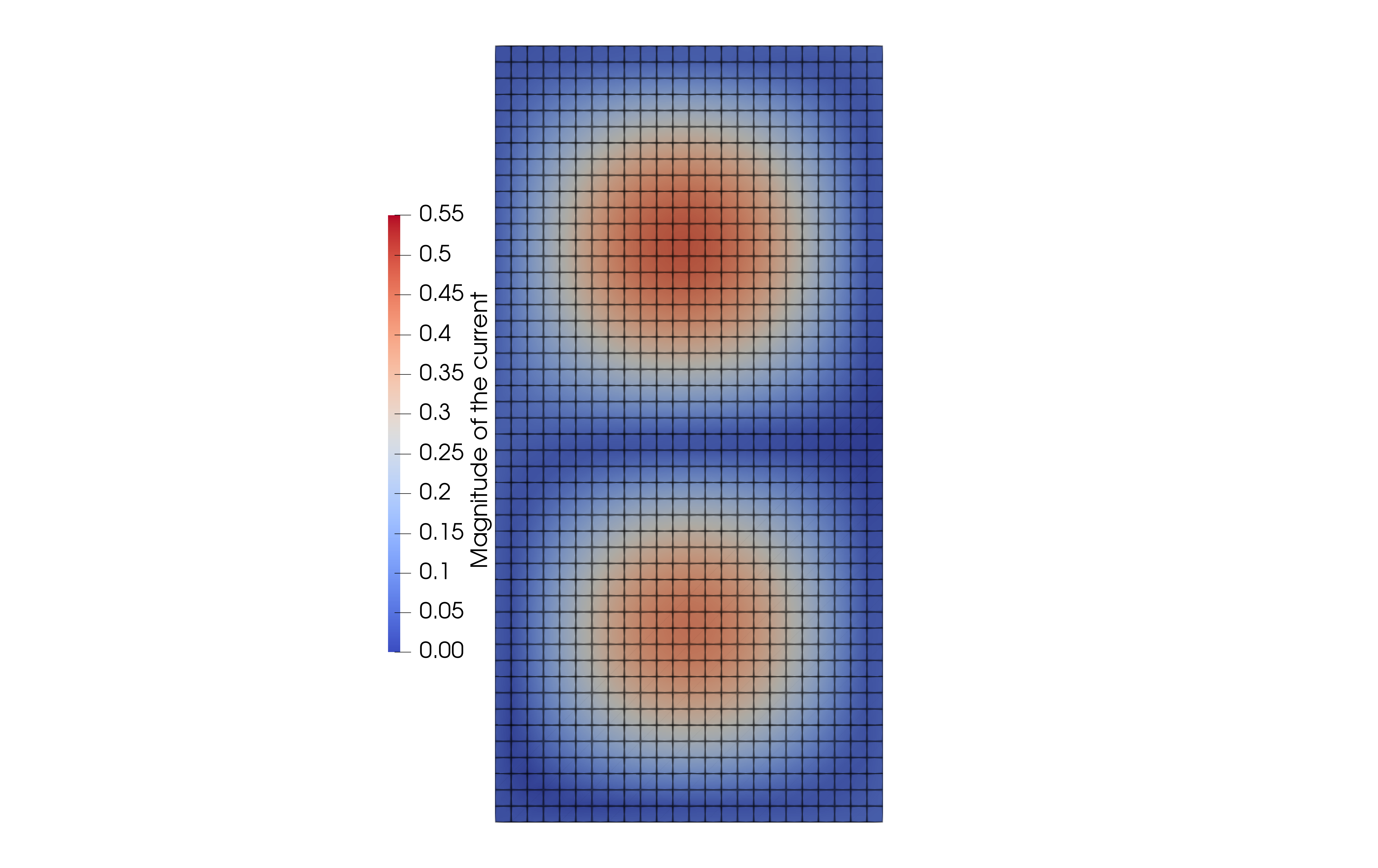}}{\qquad $\|\mathbf{j}_{\text{tar}}\|$}
\stackunder[1 pt]{\includegraphics[trim={21cm 2cm 21cm 2cm},clip,width = 1.35in]{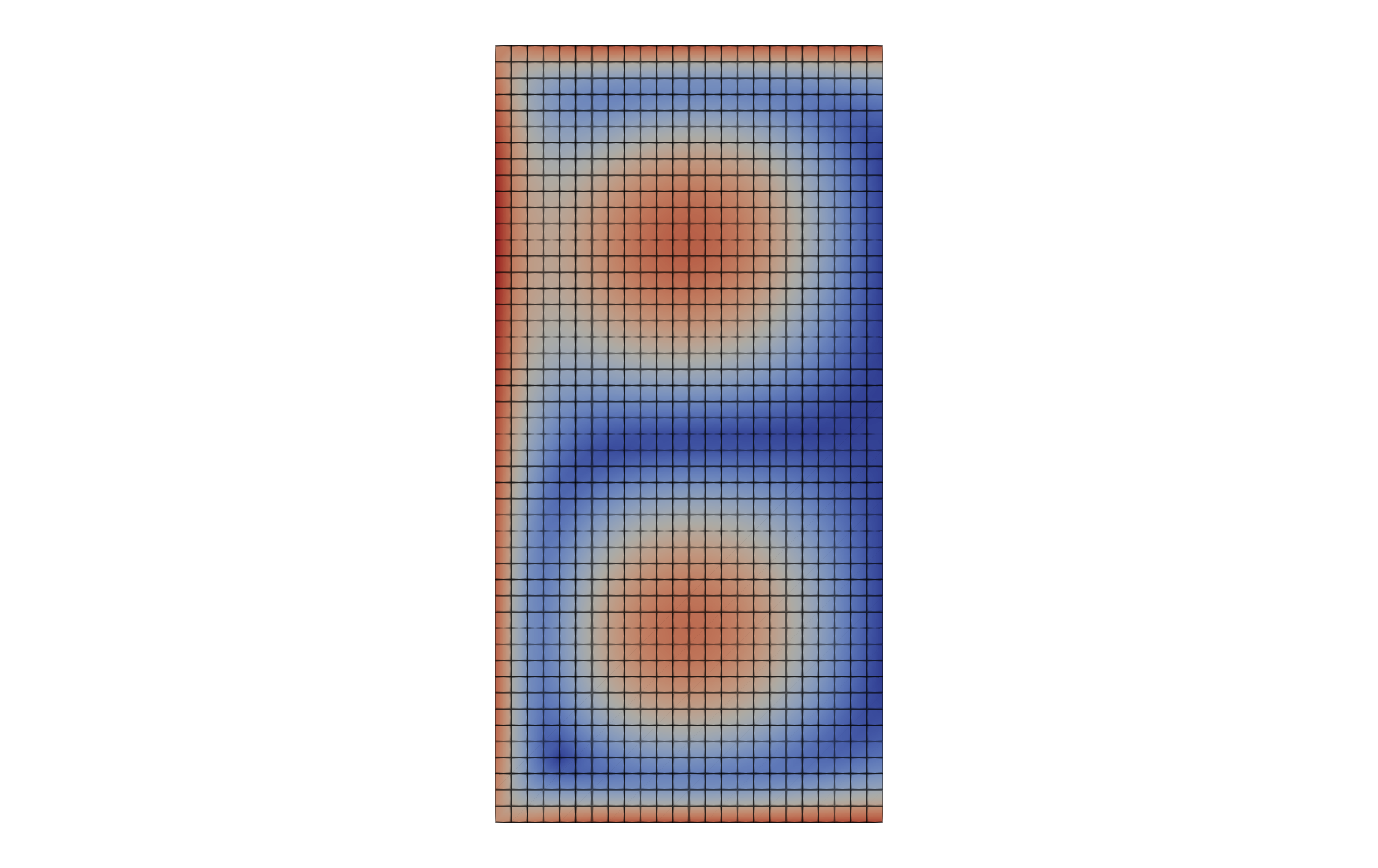}}{$\|\mathbf{j}_{\text{NS-Adj-GD}}\|$}
\stackunder[1 pt]{\includegraphics[trim={21cm 2cm 21cm 2cm},clip,width = 1.35in]{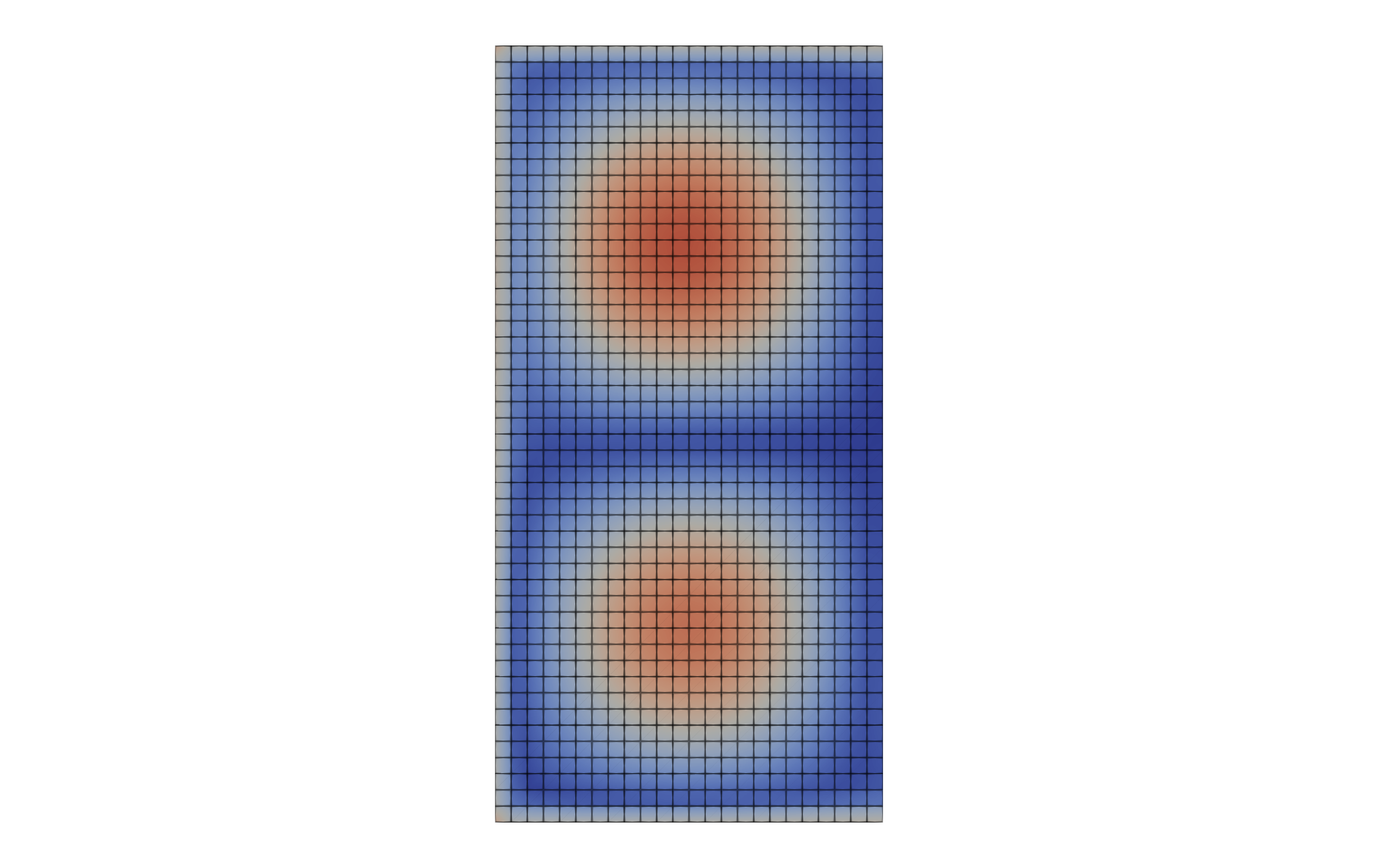}}{$\|\mathbf{j}_{\text{NS-Adj-NAG}}\|$}
\stackunder[1 pt]{\includegraphics[trim={21cm 2cm 21cm 2cm},clip,width = 1.35in]{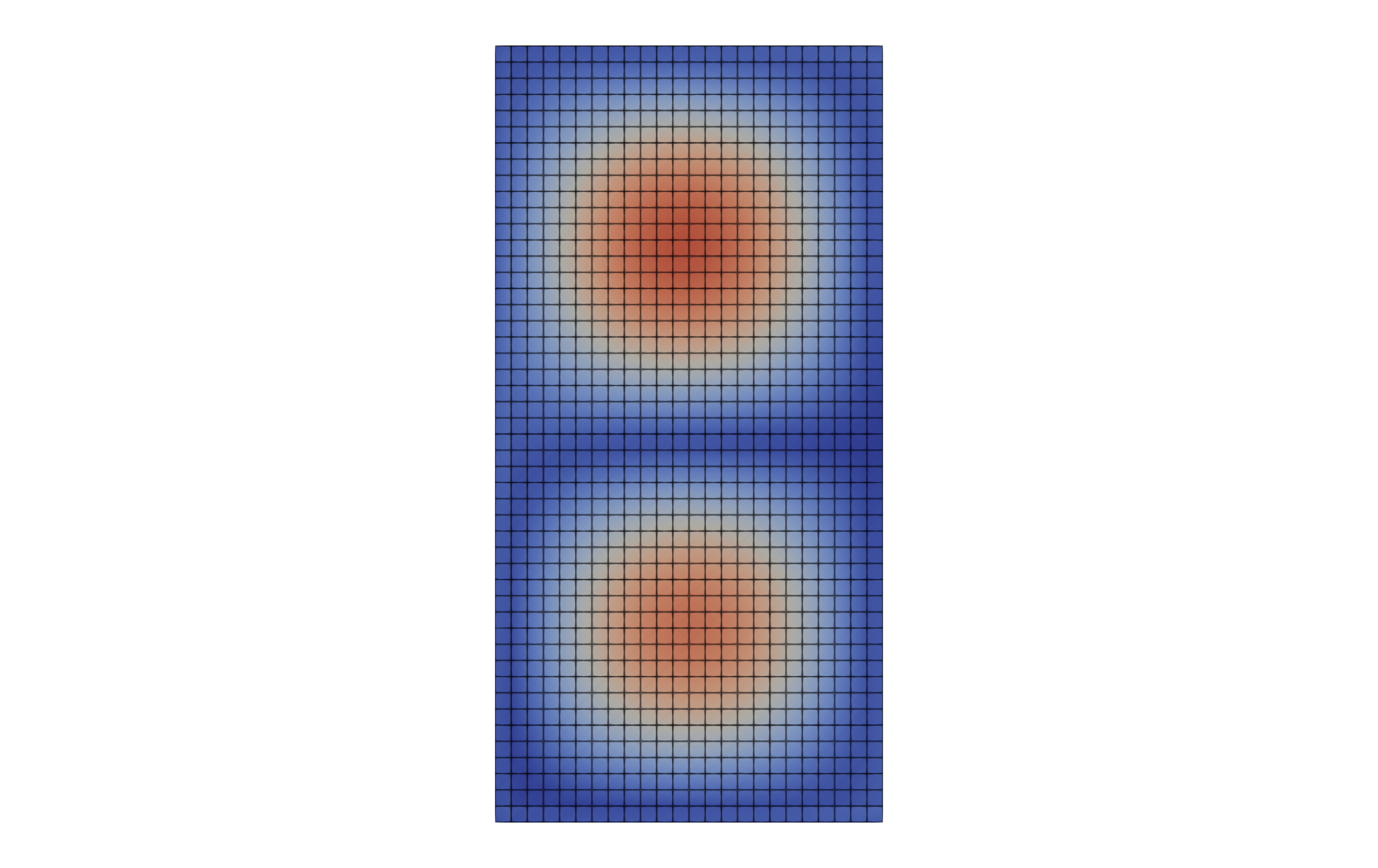}}{$\|\mathbf{j}_{\text{NS-Adj-AA}}\|$}
\caption{Comparison of the magnitude of the current of the target state ($\|\mathbf{j}_{\text{tar}}\|$) versus states optimized using NS-Adj-GD ($\|\mathbf{j}_{\text{NS-Adj-GD}}\|$), NS-Adj-NAG ($\|\mathbf{j}_{\text{NS-Adj-NAG}}\|$), and NS-Adj-AA ($\|\mathbf{j}_{\text{NS-Adj-AA}}\|$). All plots use the same scale.}\label{fig:J-solution-comp}
\end{figure}

A comparison of the scalar flux $\phi$ of the target state versus the states obtained by NS-Adj-GD, NS-Adj-NAG, and NS-Adj-AA (using the lowest values of the objective function in \Cref{fig:npg-obj}) is shown in \Cref{fig:phi-solution-comp}. Similarly, a comparison for the magnitude of the current, $\mathbf{j}$, defined as the first angular moment of $\psi$, i.e.,
$$ \mathbf{j}(x)   \coloneqq  \int_{\mathbb{S}^2} \bsO \psi(x,\bsO) d\Omega, $$
is shown in \Cref{fig:J-solution-comp}.

\subsubsection{Splitting of the initial condition for an initial value problem}\label{sec:time-int-example}
For our final example, we will consider an equality-constrained optimization problem, where the equality constraint is given by a simple linear ODE (more precisely, a time integration method applied to the ODE; for simplicity, we will use backward Euler). We will show how the nonlinear splitting framework can be used to introduce momentum into the initial condition for an optimization problem constrained by an initial value problem, and demonstrate this for a simple linear ODE.

Consider the linear initial value problem
\begin{align}\label{eq:linear-ODE}
    \dot{\bx}(t) &= A \bx(t) + \boldf(t),\ t \in [0,T], \\
    \bx(0) &= \bx_0, \nonumber
\end{align}
on a finite-dimensional vector space $X \ni \bx(t)$, $A \in L(X,X)$ is a linear map from $X$ to $X$, $\boldf: [0,T] \rightarrow X$ is a time-dependent external force, and the initial condition $\bx_0 \in X$ is given. As before, we assume that $X$ is equipped with an inner product $\langle\cdot,\cdot\rangle$ and denote the adjoint of the operator $A$ with respect to this inner product by $A^*$. Consider the backward Euler method applied to \eqref{eq:linear-ODE}, using a stepsize $\Delta t = T/N$,
\begin{align}\label{eq:backward-Euler}
    \bx_{n+1} - \bx_n &= \Delta t A \bx_{n+1} + \Delta t \boldf_{n+1},\ n = 0,\dots,N-1.
\end{align}
The backward Euler method \eqref{eq:backward-Euler}, over all unknowns $\vec{\bx}:=\{\bx_1,\dots,\bx_N\}$, forms the equality constraint under consideration. We will consider an external force optimization problem
\begin{equation}\label{eq:opt-problem-external-force}
    \min_{\{\boldf_n\}_{n=1}^N} J(\bx_N) := \|\bx_N - \bx^*\|^2
\end{equation}
subject to the equality constraint given by the backward Euler method \eqref{eq:backward-Euler}. The corresponding adjoint equation is
\begin{subequations}\label{eq:adjoint-backward-Euler}
\begin{align}
    \boldsymbol{\lambda}_N &= \frac{\partial J(\bx_N)}{\partial \bx}, \\
    \boldsymbol{\lambda}_n - \boldsymbol{\lambda}_{n+1} &= \Delta t A^* \boldsymbol{\lambda}_n,\ n = N-1, \dots, 0,
\end{align}
\end{subequations}
with unknowns $\vec{\boldsymbol{\lambda}}:=\{\boldsymbol{\lambda}_0, \dots, \boldsymbol{\lambda}_{N-1}\}$.
Now, consider the following family of nonlinear splittings of the equality constraint
\begin{equation}\label{eq:backward-Euler-global-splitting}
    (I-\Delta t A) \bx^{k+1}_{n+1} = (1-\alpha) \bx^{k+1}_n  + \alpha \bx^k_n + \Delta t \boldf_{n+1}^k,\ n=0,\dots,N-1.
\end{equation}
where the family is parametrized by $\alpha \in \mathbb{R}$. Of course, when $\alpha = 0$, this reduces to the original equality constraint, and thus, the corresponding nonlinear splitting method (\Cref{eq:ns-gd-time-update} below) reduces to standard adjoint gradient descent. On the other hand, for $\alpha \neq 0$, this can be interpreted as using momentum from the previous optimization iteration to modify the initial condition for each time step. The connection to momentum can be seen by rewriting $(1-\alpha) \bx^{k+1}_n  + \alpha \bx^k_n = \bx^{k+1}_n - \alpha (\bx^{k+1}_n - \bx^k_n)$, which is the standard initial condition at each timestep $\bx^{k+1}_n$ minus a momentum correction proportional to $\alpha$. As we will see in the simple example below, the momentum can stabilize the gradient descent. Interestingly, for $\alpha=1$, \eqref{eq:backward-Euler-global-splitting} can be computed completely parallel across all timesteps.

For simplicity, we will consider here just adjoint gradient descent for the problem \eqref{eq:opt-problem-external-force}. The nonlinear splitting adjoint gradient descent for a given $\alpha$, NS-Adj-GD($\alpha$), with this splitting takes the form
\begin{subequations}\label{eq:ns-gd-time-update}
\begin{align}
        (I-\Delta t A) \bx^{k+1}_{n+1} &= (1-\alpha) \bx^{k+1}_n  + \alpha \bx^k_n + \Delta t \boldf_{n+1}^k,\ n=0,\dots,N-1, \bx^{k+1}_0 = \bx_0, \\
    \boldsymbol{\lambda}_n^{k+1} - \boldsymbol{\lambda}_{n+1}^{k+1} &= \Delta t A^* \boldsymbol{\lambda}_n^{k+1},\ n = N-1, \dots, 0, \boldsymbol{\lambda}_N = \partial J(\bx^{k+1}_N)/\partial \bx \\
    \wtnabla_{\vec{\boldf}} J (\vec{\boldf}^k)[\vec{\bx}^k, \vec{\bx}^{k+1}] &= \vec{\boldsymbol{\lambda}}^{k+1}, \\
    \vec{\boldf}^{k+1} &= \vec{\boldf}^k - \gamma \wtnabla_{\vec{\boldf}} J (\vec{\boldf}^k)[\vec{\bx}^k, \vec{\bx}^{k+1}],
\end{align}
\end{subequations}
where $\vec{\boldf} := \{\boldf_1,\dots,\boldf_N\}$ denotes all of the optimization parameters. 

Consider the one-dimensional diffusion equation, where $A$ is the standard second-order central difference of the Laplacian on the domain $[0,1]$ with 50 interior nodes, with zero Dirichlet boundary conditions; the problem has an initial condition $\sin(6\pi x)$ and is run to a final time $T=1$ with $N=100$ timesteps. We optimize over the time-dependent external force $\vec{\boldf}$ initialized at zero, with target terminal state $\sin(2\pi x)$. 

\begin{figure}[!htb]
\subfloat[$\gamma = 10$]{\includegraphics[trim={3.7cm 8.5cm 4.5cm 9.15cm},clip,width = 2.7in]{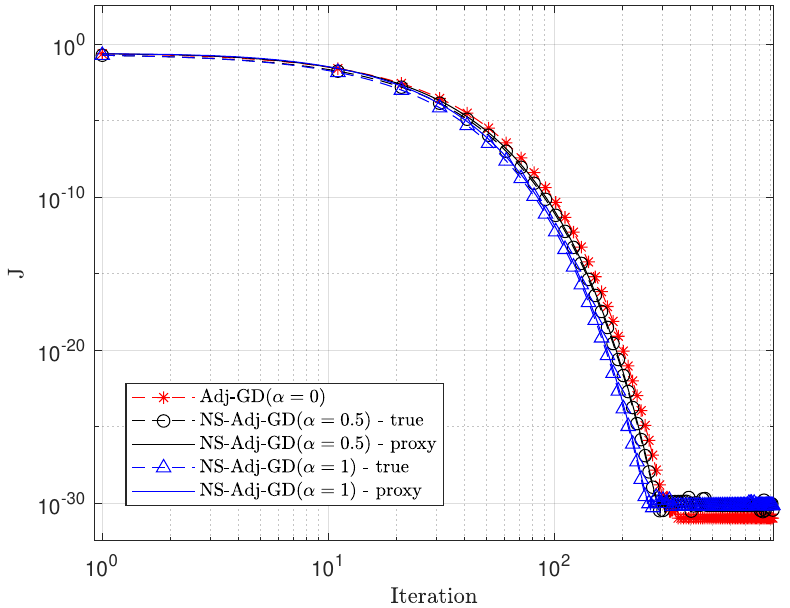}}
\subfloat[$\gamma = 50$]{\includegraphics[trim={3.7cm 8.5cm 4.5cm 9.15cm},clip,width = 2.7in]{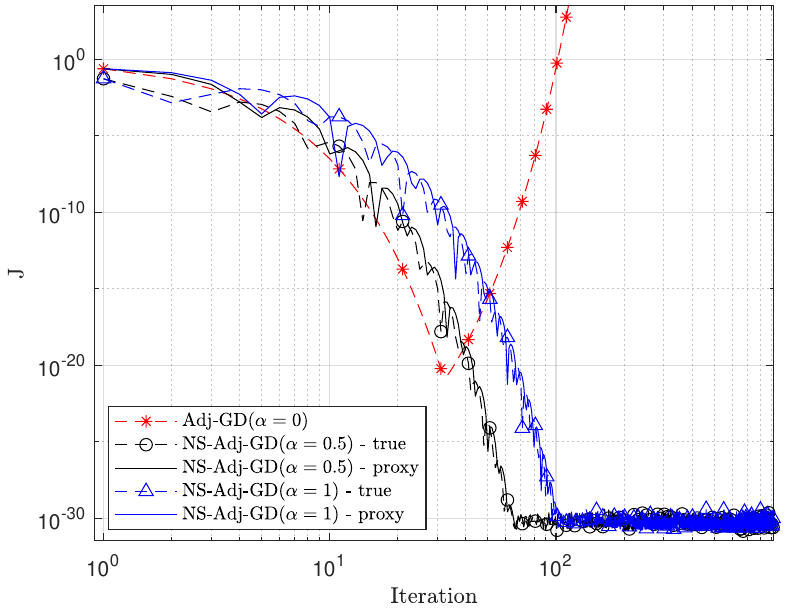}}
\caption{Comparison of NS-Adj-GD($\alpha$) for $\alpha = 0, 0.5, 1$ for various choices of the stepsize $\gamma$.}\label{fig:time-split}
\end{figure}

 In \Cref{fig:time-split}, we compare the objective versus iteration of NS-Adj-GD($\alpha$) for choices of $\alpha = 0, 0.5, 1$ (noting that $\alpha=0$ corresponds to the standard adjoint method). For the non-zero values of $\alpha$, we have included both a ``true" and a ``proxy" objective value. The proxy value at the $k^{th}$ iterate is simply $J(\bx_N^k)$; however, since we are optimizing over the external force $\vec{\boldf}$, the true value would be given by solving the equality constraint \eqref{eq:backward-Euler-global-splitting} with the current external force and then evaluating $J$. In practice, one would use the proxy as opposed to the true cost (since it does not require re-solving the equality constraint exactly), but we have included both, to show that they are in close agreement. As can be seen in \Cref{fig:time-split}, for a smaller stepsize $\gamma=10$, the three methods $\alpha=0,0.5,1$ behave similarly. On the other hand, for a larger stepsize $\gamma=50$, the standard adjoint method $\alpha=0$ becomes unstable, whereas the non-zero $\alpha$ values are able to reach the true minimum. This simple example suggests that using nonlinear splitting, e.g., in the initial conditions as considered here, can be used to stabilize an optimization method. We plan to explore this idea further in future work.

\section{Conclusion}
In this paper, we introduced a nonlinear splitting framework for gradient-based unconstrained and constrained optimization. In the unconstrained setting, we introduced a nonlinear splitting of the gradient, yielding methods with increased stability, allowing the optimization to explore the loss landscape more globally, as demonstrated by the unconstrained examples. In the constrained setting, we instead applied the nonlinear splitting to the equality constraint, which expands the space on which one can evaluate the gradient of the constrained objective, yielding a more efficient algorithm to obtain an adjoint gradient. We demonstrated the efficacy of this approach on a source optimization problem constrained by a high-dimensional kinetic equation. Interestingly, nonlinear splitting of the equality constraint can also be used to add momentum into the equality constraint; for example, we showed how the nonlinear splitting framework can be used to add momentum in the initial conditions when the constraint is an initial value problem, and provided an example of the stabilizing effect of this added momentum on a linear diffusion equation.

For future research directions, we plan to investigate the convergence and stability of nonlinear splitting optimization methods, noting that, since the nonlinear splitting framework defines a broad family of optimization methods, such analysis would require assumptions to narrow down the class of methods under consideration. It would also be interesting to investigate other applications, both unconstrained and constrained, and develop nonlinear splittings for those applications.

\section*{Acknowledgements}
BKT was supported by the Marc Kac Postdoctoral Fellowship at the Center for Nonlinear Studies at Los Alamos National Laboratory. BSS was supported by the DOE Office of Advanced Scientific Computing Research Applied Mathematics program through Contract No. 89233218CNA000001, and the National Nuclear Security Administration Interlab Laboratory Directed Research and Development program under project number 20250861ER. DBC was supported by the Center for Nonlinear Studies through the Los Alamos Summer Student Program 2025. SAS was supported by Center for Nonlinear Studies. TB was funded by the National Science Foundation under grant NSF-OIA-2327484. Los Alamos National Laboratory Report LA-UR-25-28812. 

\appendix
\section{The adjoint of the neutron transport equation}\label{sec:neutron-transport-adjoint}
We provide some additional details on the adjoint of the nonlinear neutron transport equation considered in \Cref{sec:neutron-transport}. The adjoint equation corresponding to the nonlinear neutron transport equation \eqref{eq:neutron-transport} is given by
\begin{subequations} \label{eq:neutron-transport-adjoint}
\begin{align}
    -\bsO \cdot \nabla p_\psi(x,\bsO) + \sigma_t(x) p_\psi(x,\bsO)  &- 
    \frac{\widetilde{\sigma}_s(x,\phi(x))}{4\pi}\int_{\mathbb{S}^2} p_\psi(x,\bsO) d\Omega = -g[\psi, \theta], \\
    p_\psi(x,\bsO) &= 0 \text{ for } x \in \partial D \text{ and } \bsO \cdot \mathbf{n} > 0, \label{eq:neutron-transport-adjoint-bc}
\end{align}
\end{subequations}
where $p_\psi(x,\bsO)$ is the adjoint variable, and the scattering cross section is modified for the adjoint equation to be 
\begin{align*}
     \widetilde{\sigma}_s(x,\phi(x)) &= \sigma_s(x,\phi(x)) + [D_\phi \sigma_s(x,\phi(x)) ]\phi(x) \\
     &= \sigma_s^0(x) + 3\alpha \phi(x)^2.
\end{align*}
Note that each optimization iteration involves an extra linear transport solve for the adjoint equation; these additional adjoint linear transport solves are not included in the efficiency plot, \Cref{fig:npg-obj}(B), since they are shared across all methods.

\emph{Adjoint of inflow boundary conditions.} Adjoint boundary conditions are determined by finding necessary and sufficient conditions for boundary integral terms to vanish. The adjoint is chosen to satisfy zero outflow boundary conditions \eqref{eq:neutron-transport-adjoint-bc} so no integration by parts terms arise. To see this, one linearizes the equality constraint \eqref{eq:neutron-transport} with zero inflow boundary conditions $\delta \psi (x,\bsO) = 0$ for all $x \in \partial D$ and $\bsO \cdot \mathbf{n}(x) < 0$. Given $x \in \partial D$, denote by $\mathbb{S}^2_{in}(x)  \coloneqq  \{ \Omega \in \mathbb{S}^2: \bsO \cdot \mathbf{n}(x) < 0\}$; we define $\mathbb{S}^2_{out}(x)$ analogously. Then, the boundary term between the adjoint and the linearized equation can be expressed
\begin{align*}
    \int_D \int_{\mathbb{S}^2}  & \nabla \cdot ( p_\psi(x, \bsO) \delta \psi (x,\bsO) \bsO ) d\Omega dx = \int_{\partial D} \int_{\mathbb{S}^2} p_\psi(x, \bsO) \delta \psi (x,\bsO) \bsO \cdot \mathbf{n}(x) d\Omega dS 
    \\&=  \int_{\partial D} \int_{\mathbb{S}^2_{in}(x)} p_\psi(x, \bsO) \delta \psi (x,\bsO) \bsO \cdot \mathbf{n}(x) d\Omega dS +  \int_{\partial D} \int_{\mathbb{S}^2_{out}(x)} p_\psi(x, \bsO) \delta \psi (x,\bsO) \bsO \cdot \mathbf{n}(x) d\Omega dS,
\end{align*}
where $dS$ denotes the surface measure on $\partial D$. The first term vanishes by the boundary conditions of the variational equation. The second term vanishes for arbitrary $\delta \psi$ on the outflow boundary if and only if the intensity adjoint $p_\psi$ vanishes on the outflow boundary.

\emph{Adjoint source.} The source for the adjoint equation is defined in terms of the objective, the state, and optimization parameters $\theta$. For an optimization problem
\begin{equation}
    \min_{\theta \in \Theta} J(\psi,\theta)
\end{equation}
with equality constraint given by \eqref{eq:neutron-transport}, by setting $g[\psi,\theta] = D_\psi J(\psi,\theta)$, the gradient of the objective with respect to $\nu$ can be computed via the adjoint method by
\begin{equation}
    \nabla_\theta J = \frac{\partial J}{\partial \theta} + \left( p_\psi, \frac{\partial F}{\partial \theta} \right)_{L^2_{x,\bsO}},
\end{equation}
where the equality constraint is expressed as $F = 0$ with
$$ F := \bsO \cdot \nabla \psi(x, \bsO) + \sigma_t(x) \psi(x,\bsO) - \frac{\sigma_s(x,\phi(x))}{4\pi} \int_{\mathbb{S}^2} \psi(x,\bsO) d\Omega - q(x,\bsO).$$
For the source optimization problem considered in \Cref{sec:neutron-transport},
\begin{align*}
    \min_{q(x,\bsO)} J(\psi)&=\frac{1}{2} \| \psi - \psi_{\text{tar}} \|_{L^2_{\bx,\bsO}}^2, \\
    &\text{with equality constraint } \eqref{eq:neutron-transport}, \nonumber
\end{align*}
the optimization parameter is the source $\theta = q$. The gradient of the equality-constrained objective with respect to the source is simply (minus) the adjoint variable, which can be seen by
\begin{align*}
    \nabla_q J = \left( p_\psi, \frac{\partial F}{\partial q(x,\bsO)} \right)_{L^2_{x',\bsO'}} = - p_\psi(x,\bsO),
\end{align*}
since $\partial F/\partial q(x,\bsO) = -\delta(x-x')\delta(\bsO - \bsO')$.

\bibliographystyle{plain}
\bibliography{ns_gradient.bib}

\end{document}